\documentclass[11pt]{amsart}
\usepackage{amssymb,amscd}

\newtheorem{theorem}{Theorem}[section]

\newtheorem{proposition}[theorem]{Proposition}

\newtheorem{lemma}[theorem]{Lemma}

\theoremstyle{remark}
\newtheorem{remark}[theorem]{Remark}
\newtheorem{definition}[theorem]{Definition}

\newcommand\A{\mathcal{A}}

\newcommand\be{\begin{equation}\label}
\newcommand\ee{\end{equation}}

\newcommand{\U}{\on{U}}

\newcommand{\R}{\mathbb{R}}
\newcommand{\C}{\mathbb{C}}

\newcommand{\Z}{\mathbb{Z}}

\newcommand\lie[1]{\mathfrak{#1}}
\renewcommand{\k}{\lie{k}}

\renewcommand{\b}{\lie{b}}
\newcommand{\g}{\lie{g}}

\newcommand{\n}{\lie{n}}

\renewcommand{\t}{\lie{t}}

\renewcommand{\u}{\lie{u}}

\newcommand{\on}{\operatorname}

\newcommand{\Ad}{ \on{Ad} }

\newcommand{\SU}{ \on{SU}} 
\newcommand{\SO}{ \on{SO}} 
\newcommand{\GL}{ \on{GL}}

\newcommand\dirac{/\kern-1.2ex\partial} 
\newcommand\qu{/\kern-.7ex/} 

\newcommand{\lra}{\longrightarrow}
\newcommand{\hra}{\hookrightarrow}
\newcommand{\ra}{\rightarrow}

\renewcommand{\d}{{\mbox{d}}}
\newcommand{\ol}{\overline}

\newcommand\Phinv{\Phi^{-1}}

\newcommand\sig{\sigma}

\newcommand\Om{\Omega}
\newcommand\om{\omega}

\newcommand{\f}{\frac}

\newcommand{\p}{\partial}
\renewcommand{\l}{\langle}
\renewcommand{\r}{\rangle}

\newcommand{\ti}{\tilde}

\newcommand\beqn{\begin{equation}}      
\newcommand\eeqn{\end{equation}}      
\newcommand{\ca}{\mathcal}

\newcommand{\mf}{\mathfrak}
\newcommand{\beq}{\begin{eqnarray*}}
\newcommand{\eeq}{\end{eqnarray*}}

\newcommand{\Herm}{\on{Herm}}
\newcommand{\Sym}{\on{Sym}}

\begin{document}

\title[]{Ginzburg-Weinstein via Gelfand-Zeitlin}

\author{A. Alekseev} \address{University of Geneva, Section of
  Mathematics, 2-4 rue du Li\`evre, c.p. 64, 1211 Gen\`eve 4, Switzerland}
\email{alekseev@math.unige.ch}

\author{E. Meinrenken} \address{University of Toronto, Department of
  Mathematics, 100 St George Street, Toronto, Ontario M5S3G3, Canada }
\email{mein@math.toronto.edu}

\date{\today}

\begin{abstract}
 Let $\U(n)$ be the unitary group, and $\u(n)^*$ the dual of its Lie
 algebra, equipped with the Kirillov Poisson structure. In their 1983
 paper, Guillemin-Sternberg introduced a densely defined Hamiltonian 
 action of a torus of dimension $(n-1)n/2$ on $\u(n)^*$, with moment map
 given by the Gelfand-Zeitlin coordinates. A few years later, 
 Flaschka-Ratiu described a similar, `multiplicative'  
 Gelfand-Zeitlin system for the Poisson Lie group $\on{U}(n)^*$. 
 
 By the Ginzburg-Weinstein theorem, $\on{U}(n)^*$ is isomorphic to
 $\u(n)^*$ as a Poisson manifold. Flaschka-Ratiu conjectured that one
 can choose the Ginzburg-Weinstein diffeomorphism in such a way that
 it intertwines the linear and nonlinear Gelfand-Zeitlin systems.  Our
 main result gives a proof of this conjecture, and produces a
 \emph{canonical} Ginzburg-Weinstein diffeomorphism. 
\end{abstract}

\subjclass{} 
\maketitle
\tableofcontents
\section{Introduction and statement of results}\label{sec:intro}
A theorem of Ginzburg-Weinstein \cite{gi:lp} states that for any
compact Lie group $K$ with its standard Poisson structure, the dual
Poisson Lie group $K^*$ is Poisson diffeomorphic to the dual of the
Lie algebra $\k^*$, with the Kirillov Poisson structure.  The result
of \cite{gi:lp} does not, however, give a constructive way for
obtaining such a diffeomorphism. For the case of the unitary group
$K=\U(n)$, Flaschka-Ratiu \cite{fl:mo} (see also their preprint
\cite{fl:pc1}) suggested the existence of a \emph{distinguished}
Ginzburg-Weinstein diffeomorphism, intertwining Gelfand-Zeitlin
systems on $\u(n)^*$ and $\U(n)^*$, respectively.  In this paper, we
will give a proof of the Flaschka-Ratiu conjecture.  The main result
has the following `linear algebra' implications, which may be stated
with no reference to Poisson geometry.

Let $\Sym(n)$ denote the space of real symmetric $n\times n$
matrices. For $k\le n$ let $A^{(k)}\in \Sym(k)$ denote the $k$th 
principal submatrix (upper left $k\times k$ corner) of $A\in \Sym(n)$,
and $\lambda^{(k)}_i(A)$ its ordered set of
eigenvalues, $\lambda_1^{(k)}(A)\le \cdots\le \lambda_k^{(k)}(A)$. 
The map 
\begin{equation}\label{eq:gzmap}
\lambda\colon \Sym(n)\to \R^{\f{n(n+1)}{2}},
\end{equation}
taking $A$ to the collection of numbers $\lambda_i^{(k)}(A)$ for $1\le i\le
k\le n$, is a continuous map called the \emph{Gelfand-Zeitlin map}.
Its image is the \emph{Gelfand-Zeitlin cone} $\mf{C}(n)$, cut out by the `interlacing' inequalities,
\begin{equation}\label{eq:interlacing}
\lambda_i^{(k+1)}\le \lambda_i^{(k)}\le \lambda_{i+1}^{(k+1)},\ \ 1\le
i\le k\le n-1.
\end{equation}
Now let $\Sym^+(n)\subset \Sym(n)$ denote the subset of positive definite
symmetric matrices, and define a logarithmic Gelfand-Zeitlin map 
\begin{equation}\label{eq;mu}
\mu\colon \Sym^+(n)\to \R^{\f{n(n+1)}{2}},\end{equation}
taking $A$ to the collection of numbers
$\mu^{(k)}_i(A)=\log(\lambda^{(k)}_i(A))$. Then $\mu$ is a continuous
map from $\Sym^+(n)$ onto $\mf{C}(n)$.
\begin{theorem}\label{th:real}
There is a unique continuous map $\psi\colon \Sym(n)\to \SO(n)$, with
$\psi(0)=I$, such that the map
\begin{equation}\label{gammareal}
 \gamma=\exp\circ \Ad_\psi\colon \Sym(n)\to \Sym^+(n),\ \ \
 \Ad_\psi(A)\equiv \Ad_{\psi(A)}A
\end{equation}
intertwines the Gelfand-Zeitlin maps $\lambda$ and $\mu$. In fact,
$\psi$ is smooth and $\gamma$ is a diffeomorphism. 
\end{theorem}

\emph{Remark.} For a general real semi-simple Lie group $G$ with
Cartan decomposition $G=KP$, Duistermaat \cite{du:on} proved the
existence of a smooth map $\psi\colon \mf{p}\to K$ such that 
the map $\gamma=\exp\circ \Ad_\psi\colon\mf{p}\to P$ intertwines the 
`diagonal projection' with the `Iwasawa projection'. Theorem
\ref{th:real} gives canonical maps with this property for the case
$G=\on{SL}(n,\R)$.  
\vskip.2in

\emph{Example.}  The case $n=2$ can be worked out by hand (see also
\cite[Example 3.27]{fl:mo}). Even in this case, smoothness of the map
$\gamma$ is not entirely obvious. Since $\gamma(A+tI)=e^t \gamma(A)$,
it is enough to consider trace-free matrices,
\[ A=\left(\begin{array}{cc}a &b\\ b&-a\end{array}\right).\] 
 The matrix $A$ has Gelfand-Zeitlin variables 
\[\lambda^{(2)}_1(A)=-r,\ \ \lambda^{(2)}_2(A)=r,\ \ \lambda^{(1)}_1(A)=a\] 
with $r:=\sqrt{a^2+b^2}$. Hence, the matrix $\gamma(A)$ should have 
eigenvalues $e^{-r},e^r$ and upper left entry $e^a$. This gives 
\[ \gamma(A)=\left(\begin{array}{cc} \ti{a}&\ti{b}\\ \ti{b}&\ti{c}\end{array}\right)\]
with 
\[ \ti{a}=e^a,\ \ \ti{b}=\pm \sqrt{2 e^a \cosh(r)-e^{2a}-1},\ \ \ti{c}=2\cosh(r)-e^a.\]
To obtain a continuous map, one has to take the sign of $\ti{b}$ equal
to the sign of $b$. The matrix $\psi(A)\in \SO(2)$ is a rotation
matrix by some angle $\theta(A)$. A calculation gives, 
\[ \cos(2\theta(A))=\f{a}{r}\pm \sqrt{1-\left(\f{e^a-\cosh(r)}{\sinh(r)}\right)^2}.\]
\vskip.3in

One can consider similar questions for the space $\Herm(n)$
of complex Hermitian $n\times n$-matrices, and its subset $\Herm^+(n)$
of positive definite matrices. Define surjective maps
\[ \lambda\colon \Herm(n)\to \mf{C}(n),\ \ \mu\colon \Herm^+(n)\to
\mf{C}(n)\]
in terms of eigenvalues of principal submatrices, as before. 
Let 
\[\Herm_0(n)=\lambda^{-1}(\mf{C}_0(n))\] 
denote the subset where all of the eigenvalue inequalities
\eqref{eq:interlacing} are strict. The $k$-torus 
$T(k)\subset U(k)$ of diagonal matrices acts on 
$\Herm_0(n)$ as follows,  
\begin{equation}\label{eq:torusaction}
 t\bullet A=\Ad_{U^{-1} t U}A,\ \ \ \ t\in T(k),\ A\in \Herm_0(n).
\end{equation}
Here $U\in \U(k)\subset U(n)$ is a unitary matrix such that $\Ad_{U}A^{(k)}$ is
diagonal, with entries $\lambda_1^{(k)},\ldots,\lambda_k^{(k)}$. The action
is well-defined since $U^{-1}t U$ does not depend on the choice of
$U$, and preserves the Gelfand-Zeitlin map \eqref{eq:gzmap}.  The actions
of the various $T(k)$'s commute, hence they define an action of the
\emph{Gelfand-Zeitlin torus}
\[ T(n-1)\times \cdots \times T(1)\cong \U(1)^{(n-1)n/2}.\]
Here the torus $T(n)$ is excluded, since the 
action \eqref{eq:torusaction} is trivial for $k=n$.  

Let $\Herm_0^+(n),\ \Sym_0(n)$ and $\Sym_0^+(n)$ denote the 
intersections of $\Herm_0(n)$ with $\Herm^+(n),\ \Sym(n)$ and
$\Sym^+(n)$. Thus $\Herm_0^+(n)=\mu^{-1}(\mf{C}_0(n))$.
\begin{theorem}\label{th:intertwines}
There is a unique continuous map 
\[ \gamma\colon \Herm(n)\to
  \Herm^+(n)\]
with the following three properties: 
\begin{enumerate}
\item\label{it:a} $\gamma$ intertwines the
  Gelfand-Zeitlin maps: $\mu\circ \gamma=\lambda$.
\item\label{it:b} $\gamma$ intertwines the
  Gelfand-Zeitlin torus actions on $\Herm_0(n)$ and $\Herm_0^+(n)$. 
\item\label{it:c} For any connected component $S$ of $\Sym_0(n)\subset \Herm(n)$, 
  $\gamma(S)\subset S$.
\end{enumerate}
In fact, $\gamma$ is a diffeomorphism, and has the following additional
properties: 
\begin{enumerate}\setcounter{enumi}{3}
\item\label{it:f} $\gamma$ is equivariant for the conjugation action
of $T(n)\subset \U(n)$,
\item\label{it:g} $\gamma(A+uI)=e^u \gamma(A)$.
\item\label{it:h} $\gamma(\ol{A})=\ol{\gamma(A)}$ (where the bar denotes
      complex conjugation).
\item\label{it:i} For $k\le n$, the following diagram commutes:
  \[ \begin{CD}\Herm(n) @>>> \Herm(k) @>>> \Herm(n)\\ @VV{\gamma}V
    @VV{\gamma}V @VV{\gamma}V \\ \Herm^+(n) @>>> \Herm^+(k)@>>>\Herm^+(n)
  \end{CD}\]
Here the left horizontal maps take a matrix to its $k$th principal
submatrix, while the right horizontal maps are the obvious inclusions
as the upper left corner, extended by $0$'s respectively $1$'s along
the diagonal.  
\end{enumerate}
\end{theorem}

Similar to the statement for real symmetric matrices, Theorem
\ref{th:real}, the map $\gamma$ can be written in the form
$\gamma=\exp\circ \Ad_\psi$ for a suitable map $\psi\colon \Herm(n)\to
\SU(n)$. To fix the choice of $\psi$, we have to impose an
equivariance condition under the Gelfand-Zeitlin torus action.  Given
$A\in \Herm_0(n)$, let $U_k\in U(k)$ be matrices diagonalizing
$A^{(k)}$, and $V_k\in U(k)$ matrices diagonalizing
$\gamma(A)^{(k)}=\gamma(A^{(k)})$. Then the Gelfand-Zeitlin action of
$t=(t_{n-1},\ldots,t_1)\in T(n-1)\times\cdots \times T(1)$ is given by
\[ t \bullet A=\Ad_{\chi(t,A)} A,\ \ \ t\bullet \gamma(A)=\Ad_{\ti{\chi}(t,A)} \gamma(A)\]
where 
\[ 
     \chi(t,A)=U_1^{-1}t_1 U_1\cdots U_{n-1}^{-1}t_{n-1}U_{n-1},\ \ 
\ti{\chi}(t,A)=V_1^{-1}t_1 V_1\cdots V_{n-1}^{-1}t_{n-1}V_{n-1}.
\]
Note that $\chi(t,A),\ \ti{\chi}(t,A)$ 
are independent of the choice of $U_i,V_i$. 
\begin{theorem}\label{th:phiext}
The map $\psi\colon \Sym(n)\to \SO(n),\ \psi(0)=I$ from Theorem \ref{th:real}
extends uniquely to a continuous (in fact, smooth) map $\psi\colon
\Herm(n)\to \SU(n)$ with the equivariance property 
\begin{equation}\label{eq:phiequiv}
\psi(t\bullet A)=\ti{\chi}(t,A)\psi(A)\chi(t,A)^{-1} 
\end{equation}
for all $A\in \Herm_0(n),\ t\in T(n-1)\times \cdots \times T(1)$. 
The map $\gamma$ from Theorem \ref{th:intertwines} is expressed in
terms of $\psi$ as $\gamma=\exp\circ \Ad_\psi$. 
Furthermore, 
\begin{enumerate}
\item \label{it:tha}
$\psi$ is equivariant for the conjugation action of 
$T(n)\subset \U(n)$. 
\item\label{it:thb}
$\ol{\psi(A)}=\psi(\ol{A})$.
\item\label{it:thc}
For all $k\le n$, the following diagram commutes, 
\[\begin{CD}
\Herm(k)@>>> \Herm(n)\\
@VV{\psi}V @VV{\psi}V\\
\SU(k)@>>> \SU(n)
\end{CD}
\]

\end{enumerate}
\end{theorem}
%

Note that the equivariance property \eqref{eq:phiequiv} of $\psi$ 
implies the equivariance of $\gamma$:
\[ \gamma(t\bullet A)=\exp(\Ad_{\psi(t\bullet
  A)\chi(t,A)}A)=\exp(\Ad_{\ti{\chi}(t,A) \psi(A)}A)=t\bullet
\gamma(A).
\]

Let us now place these results into the context of Poisson
geometry. Let $\mf{u}(n)$ be the Lie algebra of $\on{U}(n)$,
consisting of skew-Hermitian matrices, and identify
\[ \Herm(n)\cong \mf{u}(n)^*,\]
using the pairing $\l A,\xi\r=2\on{Im}(\on{tr}A\xi)$.  Then $\Herm(n)$
inherits a Poisson structure from the Kirillov-Poisson structure on
$\u(n)^*$. It was proved by Guillemin-Sternberg in \cite{gu:gc} that
the action of each $T(k)$ on $\Herm_0(n)$ is Hamiltonian, with moment
map the corresponding Gelfand-Zeitlin variables,
$(\lambda^{(k)}_1,\ldots,\lambda^{(k)}_k)$. On the other hand, the
unitary group $\U(n)$ carries a standard structure as a Poisson Lie
group, with dual Poisson Lie group $\U(n)^*$ the group of complex
upper triangular matrices with strictly positive diagonal entries. 
$\U(n)^*$ may be identified with $\Herm^+(n)$, by the map
taking the upper triangular matrix $X\in U(n)^*$ to the positive Hermitian matrix $(X^*X)^{1/2}\in
\Herm^+(n)$. Flaschka-Ratiu \cite{fl:pc1} proved that the $T(k)$
action on $\Herm_0^+(n)$ is Hamiltonian for the Poisson structure
induced from $\U(n)^*$, with moment map the logarithmic
Gelfand-Zeitlin variables $(\mu^{(k)}_1,\ldots,\mu^{(k)}_k)$.
\begin{theorem}\label{th:poisson}
The map $\gamma\colon \u(n)^*\to \U(n)^*$ described in 
Theorem \ref{th:intertwines} is a Poisson diffeomorphism.  
\end{theorem}
That is, for the group $K=\U(n)$ we have found a fairly explicit
description of a Ginzburg-Weinstein diffeomorphism, in Gelfand-Zeitlin
coordinates.  By contrast, no coordinate expressions are known for the
Ginzburg-Weinstein maps constructed in \cite{gi:lp,al:po,bo:st,en:co}.\\

\emph{Remark.}
A recent paper of Kostant-Wallach \cite{ko:gel} studies in detail the
holomorphic (i.e., complexified) Gelfand-Zeitlin system, for the full
space $\mf{gl}(n,\C)$. It may be interesting to consider a nonlinear
version of the holomorphic system, and to generalize our results to
that setting.\\

{\bf Acknowledgment.} We would like to thank Henrique Bursztyn for
helpful discussions.

\section{Uniqueness of the map $\gamma$}\label{sec:linalg}
In this Section, we construct a map $\gamma$ over the open dense
subset 
$\Herm_0(n)=\lambda^{-1}(\mf{C}_0(n))$, satisfying all the properties
listed in Theorem \ref{th:intertwines}. The
existence of a smooth extension to $\Herm(n)$ will be proved in the
subsequent sections.  We denote by $T_\R(k)=T(k)\cap \on{O}(k)\cong
(\Z_2)^k$ the `real part' of the torus.  The action of the
Gelfand-Zeitlin torus on $\Herm_0(n)$ restricts to an action of
$T_\R(n-1)\times\cdots T_\R(1)\cong (\Z_2)^{(n-1)n/2}$ on
$\on{Sym}^+(n)$.  The following facts concerning the Gelfand-Zeitlin map
are standard; we include the proof since we are not aware of a convenient
reference.
\begin{proposition}
  The restriction of the Gelfand-Zeitlin map to $\Herm_0(n)$ defines a
  principal bundle 
\begin{equation}\label{eq:c} \lambda\colon \Herm_0(n)\to \mf{C}_0(n)\end{equation}
  with structure group the Gelfand-Zeitlin torus $T(n-1)\times \cdots
  \times T(1)$. It further restricts to a principal bundle 
\begin{equation}\label{eq:r} \lambda\colon \Sym_0(n) \to \mf{C}_0(n)\end{equation}
with structure group $T_\R(n-1)\times\cdots \times T_\R(1)$.
Similarly for the restriction of the logarithmic Gelfand-Zeitlin map
$\mu\colon \Herm^+(n)\to \mf{C}(n)$ to $\Herm_0^+(n)$ and $\Sym_0^+(n)$.
\end{proposition}
\begin{proof}
 Consider the commutative diagram, 
\begin{equation}
\begin{CD}
\Herm_0(n) @>>> \mf{C}_0(n)\\
@VVV @VVV\\
\Herm_0(n-1) @>>> \mf{C}_0(n-1)
\end{CD} 
\end{equation}
where the horizontal maps are the Gelfand-Zeitlin maps, the left
vertical map is $A\mapsto A^{(n-1)}$, and the right vertical map
$\mf{C}_0(n)\to \mf{C}_0(n-1)$ is the obvious projection, forgetting
the variables $\lambda^{(n)}_i$. The Gelfand-Zeitlin map
$\Herm_0(n)\to \mf{C}_0(n)$ factorizes as
\begin{equation}\label{eq:iterate}
\Herm_0(n) \to \Herm_0(n-1)\times_{\mf{C}_0(n-1)} \mf{C}_0(n) \to 
\mf{C}_0(n),
\end{equation}
where the middle term is the fiber product.  By induction, we may
assume that the map $\Herm_0(n-1) \to \mf{C}_0(n-1)$, and hence the
last map in \eqref{eq:iterate}, 
is a principal bundle for the action of the Gelfand-Zeitlin torus $T(n-2)\times
\cdots \times T(1)$. 
Hence, it suffices to
show that the first map in \eqref{eq:iterate} is a principal $T(n-1)$
bundle for the Gelfand-Zeitlin action.  Thus, let $\lambda^{(k)}_i,\ 1\le i\le k\le n$ be the
components of a given point $\lambda\in\mf{C}_0(n)$, and let
$A^{(n-1)}\in\Herm_0(n-1)$ be a matrix with Gelfand-Zeitlin parameters
$\lambda^{(k)}_i$ for $1\le i\le k\le n-1$. Let us try to find
$b_1,\ldots,b_{n-1}\in\C$ and $c\in\R$ such that the matrix
\begin{equation}\label{eq:blockform}
 A=\left(\begin{array}{cc}A^{(n-1)} & b\\b^*&c\end{array}\right)
\end{equation}
has eigenvalues $\lambda^{(n)}_i$. (Here $b$ denotes the $(n-1)\times
1$-matrix with entries $b_i$.)  Choose $U\in U(n-1)$ such that the
matrix $\Lambda^{(n-1)}=U A^{(n-1)} U^{-1}$ is diagonal, with entries
$\lambda^{(n-1)}_i$ down the diagonal. Then
\[ UAU^{-1}=\left(\begin{array}{cc}\Lambda^{(n-1)} & \ti{b}\\\ti{b}^*&c\end{array}\right)
\]
where $\ti{b}=U\,b$. (As before, we think of $\U(k)$ for $k\le n$ 
as a subgroup of $\U(n)$, using the inclusion as the upper left corner.) 
The characteristic polynomial 
$\det(A-uI)$ is therefore given by, 
\[ \det(A-uI)=(c-u)
\prod_j(\lambda^{(n-1)}_j-u) -\sum_i |\ti{b}_i|^2 \prod_{j\not=i}(\lambda^{(n-1)}_j-u).\]
Setting this equal to $\det(A-uI)=\prod_r(\lambda^{(n)}_r-u)$, and
evaluating at $u=\lambda^{(n-1)}_i$ and at $u=\lambda_n^{(n)}$, one
finds 
\[ 
|\ti{b}_i|^2=-\f{\prod_r(\lambda^{(n)}_r-
  \lambda^{(n-1)}_i)}{\prod_{j\not=i}(\lambda^{(n-1)}_j-\lambda^{(n-1)}_i)},\
  \ \ c=\lambda_n^{(n)}-\sum_i
  \f{\prod_{r\not=n}(\lambda^{(n)}_r-\lambda^{(n-1)}_i)}{\prod_{j\not=i}(\lambda_j^{(n-1)}-\lambda_i^{(n-1)})}.\]
The eigenvalue inequalities ensure that the right hand side of the
expression for $|\ti{b}_i|^2$ is $>0$. This shows that the first map
in \eqref{eq:iterate} is onto.  Furthermore, since $c$ is uniquely
determined while $\ti{b}_i$ are determined up to a phase,
this map defines a principal $T(n-1)$ bundle.  Since left matrix
multiplication of $\ti{b}$ by an element of $T(n-1)$ is exactly the
Gelfand-Zeitlin action, the proof for $\Herm_0(n)$ is complete. The
proof for $\on{Sym}_0(n)$ is similar, considering only matrices with
entries in $\R$.  The parallel statements for the map $\mu$ are a
direct consequence of the statements for $\lambda$.
\end{proof}

\begin{lemma}\label{cor:ext1}\label{lem:A}
There exists a unique continuous map $\gamma\colon \Herm_0(n)\to
\Herm_0^+(n)$, satisfying \eqref{it:a}--\eqref{it:c} from Theorem 
\ref{th:intertwines}. Furthermore, this map also has the Properties
\eqref{it:f}--\eqref{it:i} from Theorem \ref{th:intertwines}.  
\end{lemma}
\begin{proof}
The choice of a connected component $S\subset \Sym_0(n)$ defines a
cross-section, hence a trivialization, of the principal 
bundle $\lambda\colon \Herm_0(n)\to \mf{C}_0(n)$. The intersection 
\[ S^+=S\cap \Sym_0^+(n)\]
is a connected component of $\Sym_0^+(n)$, which likewise trivializes
the bundle $\mu\colon\Herm_0^+(n)\to \mf{C}_0(n)$.  Thus, we obtain a
unique map $\gamma$ satisfying \eqref{it:a}--\eqref{it:b}, with
$\gamma(S)=S_+$ for the given $S$. By equivariance, the property
$\gamma(S)=S_+$ holds true for \emph{all} components $S\subset
\Sym_0(n)$, which gives \eqref{it:c}.  We claim that Property
\eqref{it:f} follows from \eqref{it:b}. Indeed, the Gelfand-Zeitlin
action of any element in $t_k\in Z(U(k))\subset T(k)\subset T(n)$
coincides with the conjugation action, since the functions
$\chi,\ti{\chi}$ in \eqref{eq:phiequiv} are simply
$\chi(t_k,A)=\ti{\chi}(t_k,A)=t_k$. The collection of these subgroups
$Z(U(k))\cong \U(1)$ of $T(n)$, together with the center $Z(\U(n))$
(which acts trivially) generate $T(n)$, proving the claim.  Properties
\eqref{it:g} and \eqref{it:h} follow from the uniqueness, since the
maps
\[A\mapsto e^{-u}\gamma(A+uI),\ \  \ \ A\to \ol{\gamma(\ol{A})}\]
satisfy \eqref{it:a}--\eqref{it:c}. Finally \eqref{it:i} holds by
the commutativity of the diagram 
  \[ \begin{CD}\Herm(n) @>>> \Herm(k) @>>> \Herm(n)\\ @VV{\lambda}V
    @VV{\lambda}V @VV{\lambda}V \\ \mf{C}(n) @>>> \mf{C}(k)@>>>\mf{C}(n)
  \end{CD}\]
and of the similar diagram for map $\mu$. 
\end{proof}

\begin{lemma}\label{lem:B}
  There is a continuous function $\psi\colon \Sym_0(n)\to \SO(n)$,
  with the property that the map $\gamma=\exp\circ \Ad_\psi\colon
  \Sym_0(n)\to \Sym_0^+(n)$ intertwines the Gelfand-Zeitlin maps. The
  map $\psi$ is uniquely determined by the additional condition
  $\psi(uA)\to I$ for $u\to 0$. 
\end{lemma}
\begin{proof}
  We have seen above that there is a unique continuous map
  $\gamma\colon \Sym_0(n)\to \Sym_0^+(n)$ which intertwines the
  Gelfand-Zeitlin maps and satisfies $\gamma(S)=S^+$ for any
  component $S\subset \on{Sym}_0(n)$.  Since $\gamma(A)$ and $\exp(A)$
  have the same eigenvalues, and since $S\cong \mf{C}_0$ is
  contractible, one can always choose a continuous map
  $\psi\colon \Sym_0(n)\to \SO(n)$ with $\gamma=\exp\circ \Ad_\psi$.

  Conversely, suppose $\psi\colon \Sym_0(n)\to \SO(n)$ is a continuous
  map, such that $\gamma=\exp\circ \Ad_\psi$ intertwines the
  Gelfand-Zeitlin maps. Suppose $\psi(uA)\to I$ for $u\to 0$.  We will
  show that (i) the map $\gamma$ has the property $\gamma(S)=S^+$
  for any connected component $S$, and (ii) the map $\psi$ with these
  properties is unique. Proof of (i): It suffices to show that the 
  restriction of $\gamma$ to $\Sym_0(n)$ is homotopic to the identity
  map of $\Sym_0(n)$. Define 
\[ [0,1]\times \Sym_0(n)\to \Sym_0(n),\ (u,A)\mapsto A_u=\begin{cases}
A& \mbox{ for }u=0\\
\f{1}{u}(\gamma(uA)-I)+uI& \mbox{ for }0<u\le 1
\end{cases}
\]
This is a well-defined continuous map since 
\[ \lim_{u\to 0}\Big( \f{1}{u}(\gamma(uA)-I)+uI\Big)=
\lim_{u\to 0}\Big(\f{1}{u}\big(\exp(\Ad_{\psi(uA)}A)-I\big)\Big)=A\]
Furthermore $A_u\in \Sym_0(n)$, since 
$\Sym_0(n)$ is invariant under
scalar multiplication by nonzero numbers, as well as under addition of
a scalar multiple of the identity matrix. Clearly $A_1=\gamma(A)$.
Proof of (ii): 
Suppose $\psi'\colon \Sym_0(n)\to \SO(n)$ is another map with
$\gamma(A)=\exp(\Ad_{\psi'(A)}A)$ and $\lim_{u\to 0}\psi'(uA)=I$. Then
$\psi'(A)=\psi(A)\chi(A)$ where $\chi(A)$ centralizes $A$ and
$\lim_{u\to 0}\chi(uA)=I$. Since the centralizer subgroup
$\on{O}(n)_A$ of any $A\in \Sym_0(n)$ is discrete, and
$\on{O}(n)_A=\on{O}(n)_{uA}$ for $u>0$, we have
$\chi(A)=\chi(uA)\xrightarrow[u\to 0]{} I$. Thus $\chi(A)=I$, proving
uniqueness of $\psi\colon \Sym_0(n)\to \SO(n)$. 
\end{proof}
Note that we have not yet shown that it is actually possible to
satisfy the normalization condition $\lim_{u\to 0}\psi(uA)=I$. This
can be proved `by hand', but will in any case be automatic for the map 
constructed below (cf. Section \ref{subsec:un}).  

\begin{lemma}
The map $\psi\colon \Sym_0(n)\to \SO(n),\ \lim_{u\to 0}\psi(uA)=I$ described in Lemma 
\ref{lem:B} admits a unique extension
$\psi\colon \Herm_0(n)\to \SU(n)$ with the equivariance property
\eqref{eq:phiequiv}.  The composition $\gamma=\exp\circ
  \Ad_\psi\colon \Herm_0(n)\to \Herm_0^+(n)$ coincides with the map
  described in Lemma \ref{cor:ext1}. Furthermore, this map also has
  the properties \eqref{it:tha} -- \eqref{it:thc} described in Theorem
  \ref{th:phiext}. 
\end{lemma}
\begin{proof}
By construction, the map 
$\gamma\colon \Sym_0(n)\to \Sym_0^+(n)$ has the equivariance 
property $\gamma(t\bullet A)=t\bullet \gamma(A)$ for all 
$t\in T_\R(n-1)\times\cdots\times T_\R(1)$. This implies the 
equivariance property \eqref{eq:phiequiv} for the map 
$\psi\colon \Sym_0(n)\to \SO(n)$, using the uniqueness part
of Lemma \ref{lem:B}.  Hence, $\psi$ admits a unique
$T(n-1)\times\cdots\times T(1)$-equivariant extension to a map
$\Herm_0(n)\to \SU(n)$, and the property $\gamma=\exp \circ \Ad_\psi$
follows by equivariance. Let us now check the additional properties
from Theorem \ref{th:phiext}.

(a) As mentioned above, the Gelfand-Zeitlin action of  
$Z(\U(k))\subset T(k)\subset T(n)$ for $k<n$ coincides with the action by conjugation. Hence, \eqref{eq:phiequiv} gives 
$\psi(\Ad_{t_k}A)=\Ad_{t_k}\psi(A)$ for $t_k\in Z(U(k))$. Since the collection of these subgroups, together with $Z(U(n))$, generate $T(n)$ it follows
that $\psi$ is $T(n)$-equivariant.

(b) $\ol{\psi(A)}=\psi(\ol{A})$
follows from the uniqueness properties, since both $\psi$ and 
$A\to \ol{\psi(\ol{A})}$ are $T(n-1)\times\cdots \times T(1)$-equivariant 
extensions of the given map over $\Sym_0(n)$.

(c) Let $\psi^{(k)}\colon \Herm(k)\to \SU(k)$ denote the analogue of
the map $\psi$, for given $k<n$. Since $\psi$ is
equivariant for the conjugation by $T(n)$, it is in particular
equivariant for the subgroup $T(n-k)$ embedded as the \emph{lower
right} corner. Since $\Herm_0(k)$ (as the upper left corner) is fixed
under this action, it follows that the restriction of $\psi$ takes
values in $\on{S}(\on{U}(k)\times T(n-k))$.  Similarly, the
restriction to $\Sym_0(k)$ takes values in $\on{S}(\on{O}(k)\times
T_\R(n-k))$. Since $T_\R(n-k)$ is discrete, the property $\lim_{u\to
0}\psi(uA)=I$ implies that $\psi|\Sym_0(k)$ must take values in
$\SO(k)$. From the uniqueness properties, it therefore follows that it
coincides with $\psi^{(k)}|\on{Sym}_0(k)$. The more general statement
$\psi|\Herm_0(k)=\psi^{(k)}|\Herm_0(k)$ now follows by equivariance
under the Gelfand-Zeitlin action of $T(k-1)\times \cdots \times T(1)$.
\end{proof}

To complete the proof of Theorems \ref{th:real}, \ref{th:intertwines},
and \ref{th:phiext}, it suffices to find a smooth map
$\psi\colon \Herm(n)\to \SU(n),\ \psi(0)=I$ with the following
Properties:
\begin{enumerate}
\item[(i)] $\psi(\ol{A})=\ol{\psi(A)}$, 
\item[(ii)] $\gamma =\exp\circ \Ad_\psi$ is a diffeomorphism intertwining the
Gelfand-Zeitlin maps and the Gelfand-Zeitlin torus actions. 
\item[(iii)] $\psi$ has the equivariance property \eqref{eq:phiequiv}. 
\end{enumerate}
The construction of a map $\psi$ with these properties, using
Poisson-geometric techniques, will be finished in Section
\ref{subsec:un}.

\section{Poisson-geometric techniques}
In this Section we discuss various tools that are needed for our construction of Ginzburg-Weinstein 
diffeomorphisms.
\subsection{Bisections}\label{subsec:groupaction}
Suppose 
\begin{equation}\label{eq:action}
 \A\colon K\to \on{Diff}(M)
\end{equation}
is an action of a Lie group $K$ on a manifold $M$. We will often write
$k.x:=\A(k)(x)$ for $k\in K,\ x\in M$. Consider the action groupoid  
\[ K\times M\rightrightarrows M\]
with face maps $\partial_0(k,x)=x$ and $\partial_1(k,x)=k.x$. 
A \emph{bisection} \cite[Chapter 15]{ca:ge} of $K\times
M\rightrightarrows M$ is a submanifold $N\subset K\times M$
such that both maps
$\partial_i$ restrict to diffeomorphisms $N\to M$. Any bisection has
the form $N=\{(x,\psi(x))|\,x\in M\}$ where $\psi\in C^\infty(M,K)$ is
a map such that 
\[ \A(\psi)(x):=\A(\psi(x))(x)\]
defines a diffeomorphism $\A(\psi)\in\on{Diff}(M)$. Henceforth, we will refer to the map $\psi$
itself as a bisection. \footnote{It can be shown that a smooth map $\psi\colon M\to
  K$ is a bisection, if and only if for all $x\in M$, the map
  $k\mapsto \psi(k.x)k$ is a diffeomorphism of $K$. In this case,
  $\psi^{-1}(x)=:h$ is obtained as the unique solution of
  $\psi(h.x)h=1$.}  Let $\Gamma(M,K)\subset C^\infty(M,K)$ denote the
set of bisections. The map 
\begin{equation}\label{eq:act1}
 \Gamma(M,K)\to \on{Diff}(M),\ \psi\mapsto \A(\psi)
\end{equation}
is a group homomorphism for the following product 
on $\Gamma(M,K)$, 
\[ (\psi_1\odot\psi_2)(x)=\psi_1(\A(\psi_2)(x))\psi_2(x).\] 
The inverse of a bisection $\psi$ for this product is given by 
\[ \psi^{-1}(x)=\psi(\A(\psi)^{-1}(x))^{-1}.\] 
The group homomorphism \eqref{eq:act1} extends the action \eqref{eq:action}
of $K$, and has kernel the bisections satisfying $\psi(x)\in K_x$ for all $x\in M$. 
%
%
%
%
For later reference we note the following easy fact: 
\begin{lemma}\label{lem:easy}
Suppose $\psi\in \Gamma(M,K)^K$ is an equivariant bisection (that is, 
$k\odot \psi=\psi\odot k$ for all $k\in K$). Then
$\psi\odot\phi=\phi\odot\psi$ for all bisections $\phi$ satisfying
$\A(\psi)^*\phi=\phi$. Furthermore,
$(\phi\odot\psi)(x)=\phi(x)\psi(x)$. 
\end{lemma}
\begin{proof}
Since $\A(\psi)^*\phi=\phi$, the product $(\phi\odot\psi)(x)$ coincides
with the pointwise product $\phi(x)\psi(x)$. On the other hand, since $\psi$
is $K$-equivariant, 
\[ (\psi\odot \phi)(x)=\psi(\A(\phi)(x))\phi(x)=\Ad_{\phi(x)}(\psi(x))\phi(x)=\phi(x)\psi(x).\]
\end{proof}

The Lie algebra $\Gamma(M,\k)$ corresponding to $\Gamma(M,K)$ may be
described as follows. 
Let
\begin{equation}\label{eq:liealg1}
 \k\to \mf{X}(M),\ \xi\mapsto \A(\xi)
\end{equation}
denote the infinitesimal generators of the $K$-action, i.e. $\A(\xi)$ is
the vector
field with flow \footnote{In this paper, we follow the convention that
the flow $F_t$ of a (possibly time dependent) vector field $X_t$ is
defined in terms of its action on functions by 
$(X_t f)(F_t^{-1}(x))=\f{\p}{\p t}f(F_t^{-1}(x))$. 
The Lie derivative $L_{X_t}$ on differential
forms is then characterized by 
$F_t^*\circ L_{X_t}=-\f{\p}{\p t}F_t^*$.}    
$F_t=\A(\exp(t\xi))$. Then \eqref{eq:liealg1}
is a homomorphism of Lie algebras.  For
$\beta\in C^\infty(M,\k)$ let $\A(\beta)\in\mf{X}(M)$ be the vector
field $\A(\beta)(x)=\A(\beta(x))(x)$.  The map $\beta\mapsto \A(\beta)$ is
a Lie algebra homomorphism for the 'action algebroid' \cite{ca:ge} Lie
bracket
\begin{equation}\label{eq:bracket}
[\beta_1,\beta_2](x)=[\beta_1(x),\beta_2(x)]+(L_{\A(\beta_1)}\beta_2)(x)-
(L_{\A(\beta_2)}\beta_1)(x).\end{equation}
on $C^\infty(M,\k)$. Let $\Gamma(M,\k)$ denote the space
$C^\infty(M,\k)$ with this Lie bracket. 

To see more clearly how $\Gamma(M,\k)$ is the infinitesimal counterpart of
$\Gamma(M,K)$, it is useful to realize $\Gamma(M,K)$ as a group of
diffeomorphisms of $K\times M$. Define two commuting actions on $K\times M$, by 
setting 
\[ \ti{\A}(k)(h,x)=(hk^{-1},k.x),\ \ \ \ti{\A}'(k)(h,x)=(kh,x).\]
%
Then the map 
\[ \Gamma(M,K)\hra \on{Diff}(K\times M),\ \ 
\psi\mapsto \ti{\A}(\partial_0^*\psi)\]
identifies $\Gamma(M,K)$ with the
group of diffeomorphism of $K\times M$ which commute with the action $\ti{\A}'$ and 
preserve the $\ti{\A}$-orbits. Similarly, 
\[ \Gamma(M,\k)\hra \mf{X}(K\times M),\ \  \beta\mapsto
\ti{\A}(\partial_0^*\beta)\] 
identifies $\Gamma(M,\k)$ with the Lie algebra of vector fields on
$K\times M$ which are invariant under the action $\ti{\A}'$ and are
tangent to the $\ti{\A}$-orbits.

Let us now assume for simplicity that $K$ is compact. For any
$\beta\in \Gamma(M,\k)$, the vector field $\ti{\A}(\partial_0^*\beta)$ is
complete, since it is tangent to orbits. Hence, its time one flow exists,
and defines an element of $\Gamma(M,K)$. We have thus extended the
exponential map $\exp\colon\k\to K$ to a map 
\[ \exp\colon \Gamma(M,\k)\to \Gamma(M,K),\]
where $\psi=\exp(\beta)$ is the unique element such that 
$\ti{\A}(\partial_0^*\psi)$ is the time one flow of $\ti{\A}(\partial_0^*\beta)$. 

More generally, one can `integrate' families of maps $\beta_t\in
\Gamma(M,\k)$ to families of bisections $\psi_t$, by viewing $\beta_t$
as a time dependent vector field $\ti{\A}(\partial_0^*\beta_t)$ on
$K\times M$ and identifying $\psi_t$ with the corresponding flow
$\ti{F}_t$ on $K\times M$.  Equivalently, let $F_t$ be
the flow of the vector field $\A(\beta_t)$ on $M$. Then $F_t=\A(\psi_t)$, where the bisection
$\psi_t\in \Gamma(M,K)$ is the solution of the ordinary differential
equation on $K$,

\begin{equation}
\label{eq:betat}
\beta_t(F_t(x))=\f{\p \psi_t(x)}{\p t}\psi_t(x)^{-1},\ \ \ \ \psi_0(x)=1.
\end{equation} 

\subsection{Gauge transformations of Poisson structures}
Let $M$ be a Poisson manifold, with Poisson bivector field $\pi$.  The
group of Poisson diffeomorphism of $(M,\pi)$ will be denoted
$\on{Diff}_\pi(M)$, and the group of Poisson vector fields by
$\mf{X}_\pi(M)$. 
Let $\sig\in\Om^2(M)$ be a closed 2-form, with the
property that the bundle map
\begin{equation}
 1+\sig^\flat \circ \pi^\sharp:\ T^*M\to T^*M
\end{equation}
is invertible everywhere on $M$. (Here $\sig^\flat\colon TM\to T^*M$
and $\pi^\sharp\colon T^*M\to TM$ are the bundle maps defined by a
2-form $\sig$ and bivector field $\pi$, respectively.)  Then the formula
\begin{equation} \pi_\sig^\sharp:=\pi^\sharp\circ (1+\sig^\flat\circ \pi^\sharp)^{-1}\end{equation}
defines a new Poisson structure $\pi_\sig$ on $M$, called the
\emph{gauge transformation of $\pi$ by $\sig$} \cite{sev:po,bu:ga}.
The symplectic leaves of $\pi_\sig$ coincide with those of $\pi$,
while the leafwise symplectic forms change by the pull-back of
$\sig$.

Gauge transformations of Poisson structures arise in the context of
Hamiltonian group actions. A Poisson action $\A\colon K\to
\on{Diff}_\pi(M)$ is called \emph{Hamiltonian}, if
there exists a \emph{moment map} $\Phi\colon M\to
\k^*$, equivariant relative to the coadjoint action on $\k^*$, 
such that the generating vector fields for the action are 
\begin{equation}\label{eq:mom1} \A(\xi)=-\pi^\sharp \l\d\Phi,\xi\r.\end{equation}
The moment map condition \eqref{eq:mom1} shows in particular that 
Hamiltonian actions preserve the symplectic
leaves. From the equivariance condition, it follows that $\Phi$ is a Poisson
map.  Conversely, any Poisson map $\Phi\colon M\to \k^*$ defines a Lie
algebra action by Equation \eqref{eq:mom1}. If $K$ is connected, and
if 
the infinitesimal $\k$-action integrates to a $K$-action, then 
the latter is Hamiltonian with $\Phi$ as its moment map.

\begin{proposition} \label{prop:bisection}
Let $(M,\pi)$ be a Hamiltonian $K$-manifold with moment map $\Phi$. 
For any bisection $\psi\in \Gamma(M,K)$ let
\begin{equation}\label{eq:sigpsi}
 \sig_\psi=-\d \l\Phi,(\psi^{-1})^*\theta^L\r,\end{equation}
where $\theta^L\in \Om^1(K)\otimes\k$ denotes the left-invariant
Maurer-Cartan form.  Then $\sig_\psi$ defines a gauge transformation
of $\pi$, and 
\[\A(\psi)_*\pi=\pi_{\sig_\psi}.\]
\end{proposition}
\begin{proof}
  Since it suffices to prove this identity leafwise, we may assume
  that $\pi$ is the inverse of a symplectic form $\om$. The moment map
  condition \eqref{eq:mom1} translates into
  $\iota(\A(\xi))\om+\l\d\Phi,\xi\r=0$. We will show
  \begin{equation}\label{eq:gaugetr} \A(\psi^{-1})^* \om=\om
    +\sig_\psi, \end{equation}
  thus in particular $\om+\sig_\psi$ is symplectic.  One easily checks
  that the pull-back of $\om$ under the map $\partial_1\colon K\times
  M\to M,\ (k,x)\mapsto k.x$ is
\[ \partial_1^*\om=\om-\d\l\Phi,\theta^L\r\in\Om^2(K\times M).\]
Equation \eqref{eq:gaugetr} follows, since $\A(\psi^{-1})$ is a
composition of $\partial_1$ with the inclusion $M\to K\times M,\ 
x\mapsto (\psi^{-1}(x),x)$.
\end{proof}
We collect some other formulas for the 2-form $\sig_\psi$ which 
will become useful later. 
\begin{proposition}\label{prop:altform}
Let $(M,\pi,\Phi)$ be as in Proposition \ref{prop:bisection}. 
\begin{enumerate}
\item
For any bisection $\psi\in \Gamma(M,K)$, 
\[
\A(\psi)^*\sig_\psi=\d \l\Phi,\psi^*\theta^L\r.
\]
\item
If $\psi_1,\psi_2\in \Gamma(M,K)$ are bisections, 
\[
\sig_{\psi_1\odot\psi_2}=\sig_{\psi_1}+\A(\psi_1^{-1})^*\sig_{\psi_2}.
\]
\end{enumerate}
\end{proposition}
\begin{proof}
  (a) Using the equivariance of $\Phi$ we have
\[ \A(\psi^{-1})^*\l\Phi,\psi^*\theta^L\r=\l \A(\psi^{-1})^*\Phi,
\A(\psi^{-1})^*\psi^*\theta^L\big)\r= \l \Phi,
\Ad_{\psi^{-1}}^{-1}\big( \A(\psi^{-1})^*\psi^*\theta^L\big)\r.\]
But
$\A(\psi^{-1})^*\psi^*\theta^L=-(\psi^{-1})^*\theta^R=-\Ad_{\psi^{-1}}\big((\psi^{-1})^*\theta^L\big)$.\\
(b) From the definition $\psi_1\odot \psi_2=(\A(\psi_2)^*\psi_1) \psi_2$
we obtain, 
\[
(\psi_1\odot\psi_2)^*\theta^L=\psi_2^*\theta^L+\Ad_{\psi_2(\cdot)^{-1}}(\A(\psi_2)^*\psi_1^*\theta^L)\]
where $\psi_2(\cdot)^{-1}$ denotes the function $x\mapsto
\psi_2(x)^{-1}$. (Not to be confused with $\psi_2^{-1}(x)$.)
Therefore, 
\[ \d \l \Phi,(\psi_1\odot\psi_2)^*\theta^L\r=\d\l\Phi,\psi_2^*\theta^L\r+\A(\psi_2)^*\d \l\Phi,\psi_1^*\theta^L\r.
\]
Now apply
$\A((\psi_1\odot\psi_2)^{-1})^*=\A(\psi_1^{-1})^*\A(\psi_2^{-1})^*$
to this result, and use (a). 
\end{proof}
Proposition \ref{prop:altform}(b) shows in particular that 
\[ \Gamma_0(M,K)=\{\psi\in\Gamma(M,K)|\  \sig_\psi=0\}\]
is a subgroup of the group of bisections. By Proposition \ref{prop:bisection}, the homomorphism $\Gamma(M,K)\to
\on{Diff}(M),\ \psi\mapsto \A(\psi)$ restricts to a group
homomorphism, 
\[\Gamma_0(M,K)\to \on{Diff}_\pi(M).\]

\subsection{Moser's method for Poisson manifolds}\label{subsec:moser}
Let $(M,\pi)$ be a Poisson manifold, and $\sig_t$ a smooth family of
closed 2-forms on $M$, with $\sig_0=0$, such that
$1+\sig_t^\flat \circ \pi^\sharp$ is invertible for all $t$.  Consider
the family of gauge-transformed Poisson structures,
$\pi_t=\pi_{\sig_t}$. Suppose
\begin{equation}\label{eq:propa}
 \f{\p }{\p t}\sig_t=-\d a_t
\end{equation}
for a smooth family of 1-forms $a_t\in\Om^1(M)$, and define a time
dependent \emph{Moser vector field} $v_t\in\mf{X}(M)$ by $v_t=-\pi_t^\sharp(a_t)$.  
Assume that the time dependent vector field $v_t$ is complete
(this is automatic if the symplectic leaves of $M$ are compact), 
and let $F_t$ be the flow with initial condition $F_0=\on{id}$.
One has \cite{al:ka},
\[ \pi_t=(F_t)_*\pi.\] 
The following alternative expression for the Moser vector field is useful:
\begin{lemma} 
The Moser vector field is given by $v_t=-\pi^\sharp(b_t)$ where 
$b_t=a_t+\iota(v_t)\sig_t$. The 1-form $b_t$ satisfies, 
\begin{equation}\label{eq:propb}
 \f{\p }{\p t}\big(F_t^*\sig_t\big)=-\d(F_t^*b_t)\end{equation}
where $F_t$ is the flow of $v_t$. 
\end{lemma} 
\begin{proof}
By definition $v_t=-\pi^\sharp (\ti{b}_t)$ where $\ti{b}_t=(1+\sig_t^\flat\circ \pi^\sharp)^{-1} a_t$. 
The calculation
\[ a_t=(1+\sig_t^\flat\circ \pi^\sharp)\ti{b}_t=\ti{b}_t-\sig_t^\flat v_t
=\ti{b}_t-\iota(v_t)\sig_t
\]
shows $\ti{b}_t=b_t$. From the definition of $b_t$ and the formula
for $\d a_t$, we find
\[ \d b_t=\d a_t+L(v_t)\sig_t=-\Big(\f{\p}{\p t}-L(v_t)\Big)\sig_t
=-(F_t^{-1})^*\f{\p}{\p t}(F_t^*\sig_t).\]
\end{proof}

We will refer to $b_t$ as the \emph{Moser 1-form}. Note that for any given 
Poisson manifold $(M,\pi)$ the list of data $v_t,F_t,a_t,b_t,\sig_t,\pi_t$ is determined by 
$b_t$ (and also by $a_t$).

The following Proposition describes a situation where the twist flows 
$\A(\psi_t)$ from Section \ref{subsec:groupaction} can be viewed as Moser flows.
\begin{proposition}\label{prop:moserflow}
Suppose $K$ is a compact Lie group, and $(M,\pi)$ is a Hamiltonian Poisson $K$-manifold with moment map
$\Phi\colon M\to \k^*$. Let $\beta_t\in \Gamma(M,\k)$ define
(cf. \eqref{eq:betat}) the family of 
bisections $\psi_t\in \Gamma(M,K)$ with $\psi_0=1$. 
Then the 2-form $\sig_t$ 
determined by the Moser 1-form  
\[ b_t=\l \d \Phi,\beta_t\r\]
coincides with the form $\sig_{\psi_t}$. Hence, the Moser flow $F_t$ 
coincides with $\A(\psi_t)$, and the gauge transformed Poisson structure 
$\pi_t=\pi_{\sig_t}$ equals $\A(\psi_t)_*\pi$. 
\end{proposition}
\begin{proof}
By the moment map property, $v_t=-\pi^\sharp(b_t)=\A(\beta_t)$, with flow 
$F_t=\A(\psi_t)$. We have to verify Equation \eqref{eq:propb}.
Observe (cf. \eqref{eq:betat})
\[\d F_t^*\beta_t= 
\d\big(\f{\p \psi_t}{\p t}
\psi_t(\cdot)^{-1}\big)=\Ad_{\psi_t}\f{\p}{\p t}
\left(\psi_t^*\theta^L\right).\]
Since $F_t^*\Phi=\Ad_{\psi_t(\cdot)^{-1}}^*\Phi$ by
equivariance of the moment map, this gives
\[ F_t^* \l \Phi, \d \beta_t\r=
\l\Ad_{\psi_t(\cdot)^{-1}}^* \Phi,\d F_t^*\beta\r=
\l\Phi,\f{\p}{\p t}(\psi_t^*\theta^L)\r
=\f{\p}{\p t} \l\Phi,\psi_t^*\theta^L\r.\]
Therefore, using Proposition \ref{prop:altform}(a), 
\[ F_t^* \d b_t=-\d F_t^* \l \Phi, \d \beta_t\r =-\f{\p}{\p
  t}\d \l \Phi,\psi_t^*\theta^L\r=-\f{\p}{\p t}F_t^*\sig_t.\]
\end{proof}


\subsection{Stability of Poisson actions}

A well-known argument due to Palais 
shows that actions of compact Lie groups $K$ on compact manifolds
$M$ are stable. That is, any deformation of such an action is obtained
via conjugation by a family of diffeomorphisms of $M$. 
This result extends to the Poisson category: 
\begin{proposition}[Stability of Poisson actions of compact Lie groups]
\label{prop:action}
  Let $(M,\pi)$ be a Poisson manifold, $K$ a compact Lie group, and
  $\A_t\colon K\to \on{Diff}_\pi(M)$ a family  of $K$-actions by Poisson
  diffeomorphism of $M$. Let $w_t\in\mf{X}(M)$ be the time dependent vector field,
  given in terms of its action on functions by
\begin{equation}\label{eq:vector}
w_t=-\int_K \d k\, \A_t(k^{-1})^*\f{\p}{\p t}\A_t(k)^*
\end{equation}
where $\d k$ denote the normalized Haar measure on $K$. Then $w_t$ is a 
Poisson vector field. If the flow $F_t\in \on{Diff}_\pi(M)$ of $w_t$ exists 
(e.g. if $M$ is compact, or if the $K$-orbits are independent of $t$), 
then 
\begin{equation}\label{eq:integrated}
 \A_t(k)\circ F_t=F_t\circ \A_0(k),\ k\in K.
\end{equation}
\end{proposition}
\begin{proof}
The vector field $w_t$ given by \eqref{eq:vector} has the property, 
\[ \f{\p}{\p t}\A_t(k)^*+w_t\circ \A_t(k)^*-\A_t(k)^*\circ w_t=0.\]
Assuming that the flow $F_t$ of $w_t$ is defined, this integrates to,
\[ F_t^*\circ \A_t(k)^*\circ (F_t^{-1})^*=\A_0(k)^*.\]
which gives \eqref{eq:integrated}. Since $\A_t(k)$ are Poisson diffeomorphisms,
each vector field $w_t(k)=-\A_t(k^{-1})^*\f{\p}{\p t}\A_t(k)^*$ is a
Poisson vector field, and hence so is the $K$-average \eqref{eq:vector}.
\end{proof}

%
%
%
\emph{Remark.}
Note that if the actions $\A_t\colon K\to \on{Diff}_\pi(M)$ commute with
another (fixed) action of a compact Lie group $H$, then the vector
field $w_t$ and hence the diffeomorphisms $F_t$ are $H$-equivariant.\\

\subsection{Poisson diffeomorphisms of $\k^*$}
Of particular importance is the case $M=\k^*$, with $\A\colon K\to
\on{Diff}_\pi(\k^*)$ the
coadjoint action. We begin with the following 
simple observation: 
\begin{lemma}\label{lem:center}
  For any compact Lie group $K$, the center of the group
  $\Gamma(\k^*,K)$ of bisections is the subgroup $\Gamma(\k^*,K)^K$ of
  equivariant bisections, and is contained in the kernel of the map
  $\Gamma(\k^*,K)\to \on{Diff}(\k^*),\ \psi\mapsto \A(\psi)$.
\end{lemma}
\begin{proof}
Suppose $\psi$ is a $K$-equivariant bisection, i.e. 
$k\odot \psi=\psi\odot k$ for all $k\in K$. Equivalently, 
$\psi(k.\mu)=\Ad_k\psi(\mu)$ for all $\mu\in \k^*,k\in
K$.  Specializing to $k\in K_\mu$, this shows that $\psi(\mu)$ is in
the centralizer of $K_\mu$. 
Since $K$ is compact, this implies
$\psi(\mu)\in K_\mu$. We have thus shown $\A(\psi)=\on{id}$ for all $\psi\in 
\Gamma(\k^*,K)^K$. Now suppose $\psi,\phi\in
\Gamma(\k^*,K)$, where $\psi$ is $K$-equivariant. Then 
\[ (\psi\odot \phi)(\mu)=\psi(\A(\phi)(\mu))\phi(\mu)=\Ad_{\phi(\mu)}(\psi(\mu))\phi(\mu)=\phi(\mu)\psi(\mu)
=(\phi\odot \psi)(\mu)\]
(for the last equality, we used that $\A(\psi)=\on{id}$).  This shows
that $\Gamma(\k^*,K)^K$ is contained in the center of
$\Gamma(\k^*,K)$. The converse is obvious, since central elements
commute in particular with elements of $K$.
\end{proof}
\emph{Remark.} A similar statement holds for invariant open subsets of
$\k^*$, with the same proof. \\

Consider $k^*$ as a Hamiltonian $K$-space, with $\Phi$ the
identity map. The subgroup $\Gamma_0(\k^*,K)$ of bisections
$\psi$ with $\sig_\psi=0$ is the group of \emph{Lagrangian
bisections}. (One can show that a bisection is Lagrangian if and only
if its graph is a Lagrangian submanifold of the symplectic groupoid
$K\times\k^*=T^*K$.) The diffeomorphism $\A(\psi)$ defined by a
Lagrangian bisection is a Poisson diffeomorphism preserving symplectic
leaves.
\begin{proposition}\label{prop:exact}
The kernel of the homomorphism 
\begin{equation}\label{eq:exact}
\Gamma_0(\k^*,K)\lra \on{Diff}_\pi(\k^*),\ \ \psi\to \A(\psi)
\end{equation}
is the group of invariant Lagrangian bisections $\Gamma_0(\k^*,K)^K$,
while its image is the normal subgroup of Poisson diffeomorphisms
preserving symplectic leaves.
\end{proposition}

\begin{proof}
By Lemma \ref{lem:center} above, $\A(\psi)=\on{id}$ for all $\psi\in
\Gamma_0(\k^*,K)^K$.  Suppose conversely that $\psi\in
\Gamma_0(\k^*,K)$ is a Lagrangian bisection with
$\A(\psi)=\on{id}$. Then each $\psi_t=r_t^*\psi$ generates the trivial
action. In particular, this is true for the constant map $\psi_0\equiv
\psi(0)$. Hence $\psi(0)$ is in the center of $K$.  Replacing $\psi$ with
$\psi'=\psi(0)^{-1}\psi$, we may assume $\psi(0)=1$. Let
$\beta_t=t^{-1}r_t^*\beta \in \Gamma(\k^*,\k)^K$ be the $\k$-valued
functions generating $\psi_t$ (cf. \eqref{eq:betat}), and $b_t=t^{-2}
r_t^* b$ the associated family of closed 1-forms. Since
$\A(\psi)=\on{id}$, the vector field $v=-\pi^\sharp(b)$ is zero. Hence
$b$ is $K$-invariant, and therefore $\psi$ is $K$-equivariant.

Let $F\in \on{Diff}_\pi(\k^*)$ be any Poisson diffeomorphism preserving
symplectic leaves. (In particular, $F(0)=0$.) We have to show
$F=\A(\psi)$ for some Lagrangian bisection $\psi$.

Suppose first that $F$ is a \emph{linear} Poisson diffeomorphism of
$\k^*$. Then $F$ is dual to a Lie algebra automorphism $f\in
\on{Aut}(\k)$. Since $F$ preserves orbits, the same is true for the 
map $f$. This implies that $f$ is an
\emph{inner} automorphism, $f=\Ad_{k}$ for some $k\in K$.  Hence
$F=\Ad^*_{k}=\A(k^{-1})$ is given by a Lagrangian bisection
$\psi\equiv k^{-1}$.

Consider now the general case. 
For $t\in\R$ let $r_t\colon \k^*\to \k^*$ denote the map
$r_t(\mu)=t\mu$.  Let $F_t=r_{t^{-1}}\circ F\circ r_t$ for $t\not=0$, so that  
the limit for $t\to 0$ is the linearization $F_0=\d_0 F$ of $F$ at the
origin. Since $(r_t)_*\pi=t\pi$, each $F_t$ is a Poisson
diffeomorphism preserving leaves. By the linear case considered above, we may assume
$F_0=\on{id}$. 

Let $v_t\in \mf{X}_\pi(M)$ be the 
time-dependent vector field, given in terms of its action on functions
by $v_t=-(F_t^{-1})^*\circ \f{\p}{\p t}F_t^*$. Write $v=v_1$. Then 
$v_t=t^{-1}(r_{t^{-1}})_*v$ for $t\not=0$. 
The vector fields $v_t$ vanish to second order at $0$, since
$F_t(0)=0$ and $\d_0 F_t\equiv \on{id}$ for all $t$. In particular, $v_0=0$.
We now use the well-known fact that a Poisson vector field on $\k^*$ is 
Hamiltonian if and only if it is tangent to the symplectic leaves
(which is automatic if $\k$ is semi-simple). This follows from the 
description of the first Poisson cohomology of $\k^*$ (see e.g. \cite{gi:lp})
\[ H^1_\pi(\k^*)\cong (\k^*)^K\otimes C^\infty(\k^*)^K.\]
Hence, we may write $v=-\pi^\sharp(b)$ for some exact 1-form $b\in
\Om^1(\k^*)$. The 1-form $b$ can be normalized by requiring that its
$K$-average be zero. (Note that exact 1-forms on $\k^*$ generate the
zero vector field if and only if they are $K$-invariant.) Letting
$b_t=t^{-2}r_t^*b$, and denoting by
$\beta_t=t^{-1}r_t^*\beta\in\Gamma(\k^*,\k)$ the corresponding
$\k$-valued functions, we get
\[ v_t=-\pi^\sharp(b_t)=\A(\beta_t).\]
Let $\psi_t\in \Gamma(\k^*,K),\ \psi_0=1$ be the family of bisections obtained by
integrating $\beta_t$ (see \eqref{eq:betat}). 
We have $\psi_t=r_t^*\psi$ with $\psi=\psi_1$.
Since the 1-forms $b_t$ are closed, the corresponding 2-forms
$\sig_t=\sig_{\psi_t}$ (cf. \eqref{eq:propb} and Proposition
\ref{prop:moserflow}) vanish. That is, the
bisections $\psi_t$ are Lagrangian.  We have $F_t=\A(\psi_t)$ by
construction, and in particular $F=\A(\psi)$.
\end{proof}

\begin{remark}
  Let $(M,\pi)$ be a Poisson manifold admitting a symplectic
  realization $S\rightrightarrows M$. 
In Bursztyn-Weinstein \cite[Section
  5]{bu:pi}, the Poisson diffeomorphisms of $M$ which are
  generated by Lagrangian bisections of $S$ are referred to as
  \emph{inner automorphisms} of $M$. Proposition \ref{prop:exact}
  characterizes the inner automorphisms for the case
  $T^*K\rightrightarrows \k^*$.
\end{remark}

%
%

\begin{proposition}\label{prop:sigma}
Suppose $\sig\in \Om^2(\k^*)$ is a closed 2-form, defining a gauge
transformation of the Kirillov-Poisson structure $\pi$ on $\k^*$. Then
there exists a bisection $\psi\in \Gamma(\k^*,K)$, $\psi(0)=1$,  such that
$\sig_\psi=\sig$. In particular, $\A(\psi)_*\pi=\pi_\sig$. 
$\psi$ is unique up to
multiplication from the right by a Lagrangian bisection. If $\sig$
is invariant under the action of $H\subset K$, the bisection $\psi$
can be taken $H$-invariant.  
\end{proposition}
\begin{proof}
The assumption on $\sig$ means that the bundle map 
$A=1+\sig^\flat\circ \pi^\sharp$ is invertible everywhere. Define a smooth 
family of closed 2-forms $\sig_t$, by letting $\sig_0=0$ and 
$\sig_t=t^{-1}r_t^*\sig$ for $t\not=0$. Introduce the corresponding 
operators 
\[ A_t=1+\sig_t^\flat\circ \pi^\sharp\]
on $T^*\k^*$, connecting $A_1=A$ with $A_0=1$.  Using $(r_t)_*
\pi=t\pi$, one finds $A_t=r_t^*\circ A\circ r_{t^{-1}}^*$ for
$t\not=0$.  Since $A$ is invertible, it follows that the operator
$A_t$ is invertible for all $t$. Hence, each $\sig_t$ defines a
gauge transformation. Now let $a_t$ be the family of 1-forms, obtained
by applying the homotopy operator to $-\f{\p}{\p t}\sig_t$, and $b_t$
the corresponding family of Moser 1-forms. By Proposition
\ref{prop:moserflow}, the bisections $\psi_t$ corresponding to 
$b_t$ satisfy $\sig_{\psi_t}=\sig_t$. Thus $\psi=\psi_1$ has the
desired property $\sig_\psi=\sig$. Uniqueness of $\psi$ up to Lagrangian bisections follows
directly from Proposition \ref{prop:altform}(b). If $\sig$ is
$H$-invariant, then the bisection $\psi$ just constructed is
$H$-invariant.    
\end{proof}

For any compact, connected Lie group $K$, we denote by $Z(K)\subset K$ the identity component of
the center, and by $K^{ss}$ its semi-simple part (commutator
subgroup). Thus $\hat{K}=K^{ss}\times Z(K)\to K$ is a finite covering
of $K$, and $\mf{k}=\k^{ss}\oplus \mf{z}(\k)$ on the level of Lie algebras. 
\begin{proposition}\label{prop:momentum}
Let $K_1,K_2$ be compact Lie groups,  
and suppose $\Phi\colon \k_2^*\to \k_1^*$ is the moment map for a
Hamiltonian action $\A\colon K_1\to \on{Diff}_\pi(\k_2^*)$. 
Suppose that the composition of $\Phi$ with the projection 
$\k_1^*\to \mf{z}(\k_1)^*$ is a linear map, $\k_2^*\to
\mf{z}(\k_1)^*$. Then there exists a Lie algebra
homomorphism $\tau\colon \k_1\to \k_2$ and a Lagrangian bisection
$\psi\in \Gamma_0(\k_2^*,K_2)$ such that $\Phi=\tau^*\circ \A(\psi)$. 
\end{proposition}
\begin{proof}
  Let us first of all observe that $\Phi(0)=0$. Indeed, for the
  $\mf{z}(\k_1)^*$-component of $\Phi$ this follows by linearity,
  while for the $(\k_1^{ss})^*$-component it follows since 
  moment maps are equivariant by definition. 

  For all $t\not=0$, the scaled Poisson homomorphism $\Phi_t=r_t^{-1}\circ
  \Phi\circ r_t$ is a moment map for the scaled action
  $k\mapsto \A_t(k)=r_t^{-1}\circ \A(k)\circ r_t$. Note that the
  $\mf{z}(\k_1)^*$-component of $\Phi_t$, and hence  
  the $Z(K_1)\subset K_1$-action, do not depend on $t$.  
  The limit $\Phi_0$ for $t\to 0$
  equals the linearization of $\Phi$ at $0$, and is a moment
  map for the linearized action $\A_0$. 
By Proposition \ref{prop:action} and the subsequent Remark, there
exists a $Z(K_1)$-equivariant Poisson diffeomorphism
  $F\in\on{Diff}_\pi(\k_2^*)$ with $\A(k)=F\circ \A_0(k)\circ F^{-1},\
  k\in K_1$.
  Since moment maps for semisimple Lie groups (in this case
  $K_1^{ss}$) are unique, and since
  the $\mf{z}(\k_1)^*$-component of $\Phi$ does not depend on $t$, 
  this implies $\Phi=\Phi_0\circ F^{-1}$. By
  Proposition \ref{prop:exact}, $F^{-1}=\A(\psi)$ for some Lagrangian
  bisection $\psi$. Since $\Phi_0$ is a linear Poisson map, it is of
  the form $\Phi_0=\tau^*$ for a Lie algebra homomorphism $\tau^*$.
\end{proof}

\section{Ginzburg-Weinstein diffeomorphisms}
The main result of this Section is Theorem \ref{th:functorial},
showing that Ginzburg-Weinstein diffeomorphisms can be arranged to 
be compatible with given Poisson Lie group homomorphisms. 
\subsection{Poisson Lie groups}
We briefly review Poisson Lie groups, referring to
\cite{dr:qu,ch:gu,lu:po} for more detailed information.  A Poisson Lie
group is a Lie group $K$, equipped with a Poisson structure $\pi^K$
for which the product map is Poisson.  The linearization of $\pi^K$ at
the group unit is a Lie algebra 1-cocycle $\delta^K\colon \k\to
\wedge^2\k$, with the property that the dual map $(\delta^K)^*$
defines a Lie bracket on $\k^*$. Conversely, if $K$ is connected, the
cobracket $\delta^K$ determines $\pi^K$. One refers to the Lie algebra
$\k$ together with $\delta^K$ as the \emph{tangent Lie bialgebra} of
the Poisson Lie group $K$. The \emph{dual Poisson Lie group} $K^*$ is
the connected, simply connected Poisson Lie group with tangent Lie
bialgebra $\k^*$. If $\pi^K=0$, the dual Poisson Lie group is simply
$\k^*$ with the Kirillov Poisson structure.

A \emph{Poisson Lie group action} of $(K,\pi^K)$ on a Poisson manifold
$(M,\pi^M)$ is a $K$-action such that the action map $\partial_1\colon
K\times M\to M,\ (k,x)\mapsto k.x$ is a Poisson map.  A
\emph{$K^*$-valued moment map} \cite{lu:mo} for such an action is a
Poisson map $\Psi\colon M\to K^*$ such that the generating vector
fields are given by
\begin{equation}\label{eq:action1} 
\A(\xi)=-(\pi^M)^\sharp \Psi^*\l\theta^R,\xi\r. 
\end{equation}
Here $\theta^R\in\Om^1(K^*)\otimes\k^*$ is the right-invariant
Maurer-Cartan form on $K^*$. Equation \eqref{eq:action1} reduces to
the usual moment map condition \eqref{eq:mom1} if $K$ carries the zero
Poisson structure. According to Lu \cite{lu:mo}, \emph{any} Poisson
map $\Psi\colon M\to K^*$ defines an infinitesimal Poisson Lie group
action via \eqref{eq:action1}. In particular, the identity map of
$K^*$ defines an infinitesimal \emph{dressing action} of $K$ on $K^*$.
In nice cases, it integrates to a global action of $K$.

Compact Lie group $K$ carry a standard Poisson structure $\pi^K$
structure, constructed by Lu and Weinstein in \cite{lu:po}.  Let
$G=K^\C$ be the complexification of $K$, viewed as a real Lie group,
and $\g=\k^\C$ its Lie algebra.  Consider the Iwasawa decomposition
\[ \g=\k\oplus
\mf{a}\oplus \mf{n},\ \ G=KA\,N\]
relative to a choice of maximal torus $T\subset K$ and fundamental
Weyl chamber. That is, $\mf{a}=\sqrt{-1}\t$ while $\n$ is the direct
sum of root spaces for the positive roots.  
Let $B(\cdot,\cdot)$ be an invariant scalar product on
$\k$, and $B^\C(\cdot,\cdot)$ its complexification.  
Then $2\on{Im} B^\C(\cdot,\cdot)$ is an invariant scalar product
on $\g$, and restricts to a non-degenerate pairing between $\k$ and
the Lie algebra $\mf{a}\oplus \mf{n}$. In this way $\k^*\cong
\mf{a}\oplus \mf{n}$ acquires a Lie algebra structure, making $\k$
into a Lie bialgebra. Thus $K$ is a Poisson Lie group, with $K^*=AN$
the dual Poisson Lie group. The action of $K$ on $G$ by left
multiplication descends to the dressing action $\A_{K^*}$ on $K^*$, viewed as a
homogeneous space $G/K$. To analyze the dressing action, it is convenient to work with
the Cartan decomposition 
\begin{equation}\label{eq:cartandec}
\g=\k\oplus \mf{p},\ \ G=K\,P.\end{equation} 
where $\mf{p}=\sqrt{-1}\k$ and $P=\exp\mf{p}$. Recall that the exponential
map $\exp\colon\g\to G$ restricts to a diffeomorphism $\mf{p}\to P$. Let
$e\colon \k^*\to K^*$ be the diffeomorphism defined by the compositions,
\[ \k^*\cong \g/\k\cong \mf{p}\stackrel{\exp}{\lra} P \cong
G/K\cong K^*.\]
Then $e$ intertwines the coadjoint action $\A_{\k^*}$ with the dressing action: 
\[e\circ \A_{\k^*}(k)=\A_{K^*}(k)\circ e.\]
\emph{Example.}
 Let $K=\U(n)$, with maximal torus $T=T(n)$ and the
  usual choice of positive roots.  Then $G=\GL(n,\C)$ (viewed as a
  real Lie group), $N$ are the upper triangular matrices with $1$'s
  down the diagonal, and $A$ are the diagonal matrices with positive
  entries. Hence $K^*=AN$ are the upper triangular matrices with
  positive diagonal entries. Furthermore, $\mf{p}=\Herm(n)$ and
  $P=\Herm^+(n)$. The isomorphism $K^*\cong P$ is explicitly
  given by $X\mapsto (X^* X)^{1/2}$, and identifies the dressing
  action with the conjugation action on $\Herm^+(n)$.


\subsection{Ginzburg-Weinstein diffeomorphisms}\label{subsec:gw}
Let $K$ be a compact Lie group with the standard Poisson structure,
and consider the map $e\colon \k^*\to K^*$ constructed above. In
\cite{al:po}, it was observed that the Poisson structure 
$\pi^{\k^*}_1=(e^{-1})_*\pi^{K^*}$
is gauge equivalent to the Kirillov-Poisson structure
$\pi^{\k^*}_0=\pi^{\k^*}$. 
\begin{theorem}\label{th:alek}
  \cite{al:po} There is a canonical $T$-invariant closed 2-form
  $\sig\in \Om^2(\k^*)$, with the property 
\[ (e^{-1})_*\pi^{K^*}=\pi^{\k^*}_{\sig}.\] 
\end{theorem}
See \cite{al:li} for an explicit description of the 2-form $\sig$. 
%
We can now state a refined version of the Ginzburg-Weinstein
theorem \cite{gi:lp}.  A similar result was obtained by
Enriquez-Etingof-Marshall in \cite{en:co}, for \emph{formal} Poisson
Lie groups.\\
\begin{theorem}[Ginzburg-Weinstein diffeomorphisms]\label{th:gw}
Let $K$ be a compact Lie group with the standard Poisson structure. 
Then there exists a bisection $\psi\in \Gamma(\k^*,K)$, with
$\psi(0)=1$, such that the map 
\[ \gamma=e\circ \A(\psi)\colon \k^*\to K^*\]
is a Poisson diffeomorphism.   The
bisection $\psi$ can be chosen to be $T$-equivariant and to take values
in the semi-simple part $K^{ss}$.
\end{theorem}
\begin{proof}
By Proposition \ref{prop:sigma}, there exists a bisection $\psi\in
\Gamma(\k^*,K),\ \psi(0)=1$ with
$\sig_\psi=\sig$. For any such bisection 
$\A(\psi)_*\pi^{\k^*}=\pi^{\k^*}_\sig=(e^{-1})_*\pi^{K^*}$.
Since $\sig$ is $T$-invariant, one can take $\psi$ to be $T$-equivariant. 

The map $\psi$ lifts to a unique map $\hat{\psi}\in
\Gamma(M,\hat{K}),\hat{\psi}(0)=1$ with values in the finite 
cover $\hat{K}=K^{ss}\times Z(K)$ of $K$. Replacing $\psi$ with the $K^{ss}$-component of
$\hat{\psi}$, we arrange that $\psi$ takes values in $K^{ss}$. 
\end{proof}
\begin{definition}
A bisection $\psi\in \Gamma(\k^*,K)$ will be called a
\emph{Ginzburg-Weinstein twist} if it has the properties $\psi(0)=1$ and
$\sig_\psi=\sig$. 
\end{definition}
By Proposition \ref{prop:altform}(b), Ginzburg-Weinstein twists are 
unique up to a Lagrangian bisection. 

Ginzburg-Weinstein twists can be
used to turn ordinary $\k^*$-valued moment maps into $K^*$-valued
moment maps, and vice versa. However, the change of the moment map
produces a twisted action. 

\begin{definition}
Given an $K\to \on{Diff}(M)$ on a manifold $M$, and a bisection $\psi\in \Gamma(M,K)$, 
the \emph{$\psi$-twisted action} of $K$ on $M$ is defined as follows, 
\begin{equation}\label{eq:twistact}
\A^\psi\colon K\to \on{Diff}(M),\ \A^\psi(k)=
\A(\psi)\circ
\A(k)\circ \A(\psi)^{-1}.\end{equation}
\end{definition}

\begin{proposition}\label{prop:moment} 
Suppose $\psi\in \Gamma(\k^*,K)$ is a Ginzburg-Weinstein twist, and
let $\gamma=e\circ \A(\psi)$. Let $(M,\pi)$ be a Poisson manifold, and
$\Phi\colon M\to \k^*$ and $\Psi\colon M\to K^*$ two Poisson maps related by 
$\Psi=\gamma\circ \Phi$. 
Then $\Phi$ is the moment map for a Hamiltonian $K$-action $\A$ if and only if
$\Psi$ is the moment map for a Poisson Lie group $K$-action $\A'$. The two
actions are related as follows,
\begin{equation}
\label{eq:actions}
\A'=\A^{\Phi^*\psi^{-1}},\ \ \A=(\A')^{(e^{-1}\circ \Psi)^*\psi}.
\end{equation}
\end{proposition}
\begin{proof}
  Suppose $\Phi$ generates a $K$-action $\A$.
  We will show that $\Psi$ generates the action
  $\A'=\A^{\psi^{-1}\circ \Phi}$. 
By \cite{al:po,al:li}, the map 
\[ e\circ \Phi\colon M\to K^*\]
is the moment map for a Poisson-Lie group action of $(K,\pi^K)$ on
$M$, where $M$ is equipped with the gauge transformed Poisson structure
$\pi_{\Phi^*\sig}$. Since $\sig=\sig_\psi$, the diffeomorphism 
$\A(\Phi^*\psi^{-1})$ takes the gauge transformed Poisson structure
$\pi_{\Phi^*\sig}$ structure back to $\pi$. Furthermore,
$\A(\Phi^*\psi^{-1})$ intertwines $\A$ with the twisted action 
$\A'=\A^{\psi^{-1}\circ \Phi}$, and takes $e\circ \Phi$ to
\[ (e\circ \Phi)\circ \A(\Phi^*\psi)=
e\circ \A(\psi)\circ \Phi= \Psi .\]
It follows that $\Psi$ is a moment map for the twisted action
$\A'=\A^{\Phi^*\psi^{-1}}$ on $(M,\pi)$. Conversely, assume that $\Psi$ generates
an actions $\A'$. Then $e^{-1}\circ \Psi\colon M\to \k^*$ is a moment
map for a Hamiltonian action on $(M,\pi_{-(e^{-1}\circ \Psi)^*\sig})$.
  Applying $\A((e^{-1}\circ \Psi)^*\psi$ to restore the Poisson
  structure $\pi$, and arguing as above, we see
  that $\Phi$ is a moment map for the action
  $\A=(\A')^{(e^{-1}\circ\Psi)^*\psi}$ 
  on $(M,\pi)$. 
(Alternatively, one can also use Lemma
    \ref{lem:twisttwist} below to argue that the two formulas
    \eqref{eq:action} are equivalent.)
\end{proof}

\subsection{Functorial properties of Ginzburg-Weinstein maps}
A homomorphism of Poisson Lie groups $K_1,K_2$ is a Lie group homomorphism 
\[ \ca{T}\colon K_1\to K_2\]
which is also a Poisson map. On the infinitesimal level, a
homomorphism of Poisson Lie groups defines a homomorphism of Lie
bialgebras, $\tau\colon \k_1\to \k_2$. That is, the dual map
$\tau^*\colon \k_2^*\to \k_1^*$ is a Lie algebra homomorphism, and in
particular exponentiates to a dual Poisson Lie group homomorphism
$\ca{T}^*\colon K_2^*\to K_1^*$.  Given a Poisson Lie group action
$\A$ of $K_2$ on a Poisson manifold $M$, with moment map $\Psi\colon
M\to K_2^*$, the composition $\ca{T}^*\circ \Psi$ is a moment map for
the $K_1$-action $\A\circ \ca{T}$.

For any compact Lie group $K$ with the standard Poisson structure, the
maximal torus $T$ with the zero Poisson structure is a Poisson-Lie
subgroup. That is, the inclusion $\ca{T}\colon T\to K$ is a
Poisson-Lie group homomorphism. 
\begin{lemma}\label{lem:abeliancase}
Suppose $\psi\in\Gamma(\k^*,K)$ is a $T$-equivariant Ginzburg-Weinstein
twist, and let $\gamma=e\circ \A(\psi)$. Then the following diagram commutes: 
\[\begin{CD}
\k^* @>>{\tau^*}> \t^*\\
@VV{\gamma}V @VV{\cong}V\\
K^* @>>{\ca{T}^*}> T^*
\end{CD}\]
\end{lemma}
\begin{proof}
  Let $\ca{T}\colon T\to K$ denote the inclusion, and consider the
  Poisson map $\Upsilon\colon \k^*\to \t^*$ given as the composition
  of the Poisson maps $\gamma\colon \k^*\to K^*$ and $\ca{T}^*\colon
  K^*\to T^*\cong \t^*$. Proposition \ref{prop:moment} shows that $\gamma$
  is a moment map for the twisted $K$-action $\A^{\psi^{-1}}$, and
    hence $\Upsilon$ is a moment map for the twisted $T$-action,
    $\A^{\psi^{-1}}\circ\ca{T}$. Since $\psi$ is $T$-equivariant, the
    twisted and untwisted $T$-actions coincide. Hence, $\Upsilon$ and
    $\tau^*$ are moment maps for the same $T$-action. It follows that
    their difference is a $K$-invariant function $\k^*\to \t^*$. It
    hence suffices to show that $\Upsilon$ and $\tau^*$ coincide on
    $\t^*=(\k^*)^T\subset \k^*$ (fixed point set for the coadjoint
    action of $T$ on $\k^*$). That is, we have to show that $\Upsilon$
    restricts to the identity map of $\t^*$. 

    Since $\psi$ is $T$-equivariant, it
    takes $\t^*=(\k^*)^T$ to $T=K^T$ (fixed point set for the conjugation
    action of $T$ on $K$).  In particular, $\A(\psi)$ acts trivially
    on $\t^*$, and hence $\gamma$ coincides with $e$ on $\t^*\subset
    \k^*$. Since $e\colon \k^*\to K^*$ restricts to the natural
    identification $\t^*\cong T^*$, we conclude that $\Upsilon$
    restricts to the identity map of $\t^*$.
\end{proof}

\begin{theorem}[Compatible Ginzburg-Weinstein maps]\label{th:functorial}
  Let $K_1,K_2$ be compact Poisson Lie groups with the standard
  Poisson structure, and $\ca{T}\colon K_1\to K_2$ a Poisson Lie group
  homomorphism. Given any Ginzburg-Weinstein twist
  $\psi_{1}\in \Gamma(\k_1^*,K_1)$, there exists a
  Ginzburg-Weinstein twist $\psi_{2}\in
  \Gamma(\k_2^*,K_2)$, for which the diagram 
\begin{equation}\label{eq:commdia}
\begin{CD}
\k_2^* @>{\tau^*}>> \k_1^*\\
@VV{\gamma_{2}}V @VV{\gamma_{1}}V\\
K_{2}^* @>>{\ca{T}^*}> K_{1}^*
\end{CD}
\end{equation}
with $\gamma_{i}=e_{i}\circ \A_i(\psi_{i})$, commutes. Here $\A_i$
denotes the coadjoint action of $K_i$ on $\k_i^*$. One can
arrange that $\psi_{2}$ takes values in the semi-simple part
$K_2^{ss}$. 
\end{theorem}
\begin{proof}
We may assume, passing to a finite cover of $K_1$ if necessary, that
the semi-simple part $K_1^{ss}$ is simply
connected. We begin by choosing an arbitrary $T_2$-equivariant Ginzburg-Weinstein twist
$\psi_2$.  We will show how to modify $\psi_2$ (possibly destroying
the $T_2$-equivariance), in such a way that the above
diagram commutes. The idea is to apply Proposition \ref{prop:momentum} 
to the Poisson map 
\[ \Upsilon=\gamma_1^{-1}\circ \ca{T}^*\circ \gamma_2\colon
\k_2^*\to \k_1^*.\] 
To apply this Proposition, we have to verify that the
$\mf{z}(\k_1)^*$-component of $\Upsilon$ is given by a linear map. 
In fact, we will show that the $\mf{z}(\k_1)^*$-components of $\Upsilon$ 
and $\tau^*$ are equal. Since $\psi_2$ is $T_2$-equivariant, 
Lemma \ref{lem:abeliancase} shows that 
$\gamma_2$ restricts to the natural identification $\t_2^*\to T_2^*$. 
Similarly, $\gamma_1$ 
restricts to the natural identification $\mf{z}(\k_1)^*\to 
Z(K_1)^*=\mf{z}(\k_1)^*$, since this is true for $e_2$, and since the
action of $K_1$ (hence of $\A(\psi_1)$) on $\mf{z}(\k_1)^*$ is
trivial. Since the 
Poisson bivector of $K_2$ vanishes exactly along $T_2$, the map $\ca{T}$
must take $Z(K_1)\subset T_1$ into $T_2$. Hence, the diagram 
\[
\begin{CD}
\k_2^* @>>> \t_2^*@>>> \mf{z}(\k_1)^*\\
@VV{\gamma_{2}}V @VV{\gamma_{2}}V@VV{\gamma_{1}}V\\
K_{2}^* @>>> T_2^*@>>> Z(K_{1})^*
\end{CD}
\]
commutes, proving the claim. 

It follows in particular that the $\mf{z}(\k_1)^*$-component of
$\Upsilon$ is a moment map for the action of $Z(K_1)\subset K_1$ via
$\ca{T}$. On the other hand, the $(\k_1^{ss})^*$-component is a moment
map for some action of $K_1^{ss}$, since $K_1^{ss}$ is simply
connected. Hence, $\Upsilon$ is the moment map for a $K_1$-action.  
By Proposition \ref{prop:momentum}, there exists a 
Lagrangian bisection $\phi\in \Gamma_0(\k_2^*,K_2),\ \phi(0)=1$, with the property 
$\Upsilon\circ \A(\phi)=\Upsilon_0=\tau^*$. That is, replacing $\psi$
with $\psi'=\psi\circ \phi$ the diagram \eqref{eq:commdia} commutes.
As in the proof of Theorem \ref{th:gw}, one can arrange that the new $\psi$
takes values in $K_2^{ss}$, without changing $\A(\psi)$. 
\end{proof}

Let us call two Ginzburg-Weinstein twists $\psi_1\in\Gamma(\k_1^*,K_1)$ and $\psi_2\in
\Gamma(\k_2^*,K_2)$ \emph{compatible} (relative to $\ca{T}\colon
K_1\to K_2$) if the corresponding Ginzburg-Weinstein diffeomorphism
$\gamma_i=e_i\circ \A(\phi_i)$ define a commutative diagram
\eqref{eq:commdia}. The compatibility condition is equivalent to a
certain equivariance condition, as the following result shows. 
\begin{theorem}\label{th:equivariance}
Suppose $\psi_1\in\Gamma(\k_1^*,K_1)$ and $\psi_2\in
\Gamma(\k_2^*,K_2)$ are compatible Ginzburg-Weinstein twists, and put
\[ \hat{\psi}_1=\ca{T}\circ \psi_{1}\circ (e_1^{-1}\circ \ca{T}^*\circ e_{2})\in\Gamma(\k_2^*,K_2).\]
Then the 'ratio' $\hat{\psi}_1^{-1}\odot \psi_2\in\Gamma(\k_2^*,K_2)$ is
$K_1$-equivariant in the sense that it $\odot$-commutes with all $\ca{T}(k)$ for all $k\in
K_1$. One has the formula, 
\begin{equation}\label{eq:anotherform}
(\hat{\psi}_1^{-1}\odot \psi_2)(\mu)
=\ca{T}\big(\psi_1(\tau^*\mu)\big)^{-1}\,\psi_2(\mu).
\end{equation}
\end{theorem}
\begin{proof}
Given \emph{arbitrary} Ginzburg-Weinstein twists $\psi_1,\psi_2$, 
consider again the moment map $\Upsilon=\gamma_{1}^{-1}\circ \ca{T}^*\circ
  \gamma_{2} \colon \k_2^*\to \k_1^*$ as in the proof of Theorem
\ref{th:functorial}. Suppose $(M,\pi)$ is a Poisson manifold, and
$\Phi\colon M\to \k_2^*$ is the moment map for a Hamiltonian action
$\A\colon K_2\to \on{Diff}_\pi(M)$. Let us compute the $K_1$-action
generated by $\Upsilon\circ \Phi$. Using Proposition \ref{prop:moment} we have, 
\[ \begin{array}{crcl}
&\Phi\colon M \to \k_2^*&\mbox{ is a moment map for }& \A\\
\Rightarrow&\gamma_2\circ \Phi\colon M\to K_2^*&\mbox{ is a moment map for }& \A^{\Phi^*\psi_{2}^{-1}}\\
\Rightarrow& \ca{T}^*\circ\gamma_{2}\circ \Phi\colon M\to K_1^*&\mbox{ is a
  moment map for }& \A^{\Phi^*\psi_{2}^{-1}}\circ \ca{T}\\
\Rightarrow &\gamma_{1}^{-1}\circ \ca{T}^*\circ
  \gamma_{2}\circ\Phi \colon M\to \k_1^*&\mbox{ is a moment map for }&(\A^{\Phi^*\psi_{2}^{-1}}\circ \ca{T})^{\psi_{1}\circ
  e_{1}^{-1}\circ (\ca{T}^*\circ \gamma_{2}\circ \Phi)}.
\end{array}\]
We may re-write the result as 
\[\begin{split}
(\A^{\Phi^*\psi_{2}^{-1}}\circ \ca{T})^{\psi_{1}\circ
  e_{1}^{-1}\circ \ca{T}^*\circ \gamma_{2}\circ \Phi}
&=(\A^{\Phi^*\psi_{2}^{-1}})^{\Phi^*(\ca{T}\circ \psi_{1}\circ
  e_{1}^{-1}\circ \ca{T}^*\circ \gamma_{2})}\circ \ca{T}\\
&= (\A^{\Phi^*\psi_{2}^{-1}})^{\Phi^*(\hat{\psi}_{1}\circ
  \A(\psi_2))}\circ \ca{T}\\
&=\A^{\Phi^*(\psi_{2}^{-1}\odot \hat{\psi}_{1})}\circ \ca{T}.
\end{split}
\]
In the last line, we have used Lemma \ref{lem:twisttwist} below to
write an iterated twist as a single twist. 

Assume now that $\psi_1,\psi_2$ are compatible. The commutativity of
the diagram \eqref{eq:commdia} means that $\Upsilon =\tau^*$.  In
particular, for any Hamiltonian $K_2$-space $(M,\pi,\Phi)$, the
twisted $K_1$-action $\A^{\Phi^*(\psi_2^{-1}\odot \hat{\psi}_1)}\circ \ca{T}$ coincides with the 
untwisted action $\A\circ \ca{T}$. By definition of the twisted
action, this is equivalent to
\begin{equation}\label{eq:good}
 \A(\Phi^*(\psi_2^{-1}\odot \hat{\psi}_1))\circ \A(\ca{T}(k))=\A(\ca{T}(k))\circ
\A(\Phi^*(\psi_2^{-1}\odot \hat{\psi}_1)),\ \ k\in K_1.\end{equation}
Apply this result to $M=K_2\times\k_2^*$, with symplectic structure
coming from the identification with $T^*K_2$, and with $\Phi\colon
(k,\mu)\mapsto \mu$ the
moment map for the $K_2$-action $\A(k)(h,\mu)=(hk^{-1},k.\mu)$. Since
the map $\Gamma(\k_2^*,K_2)\to \on{Diff}(K_2\times\k_2^*),\
\psi\mapsto \A(\Phi^*\psi)$ is 1-1,
the above equation implies $(\psi_2^{-1}\odot \hat{\psi}_1)\odot \ca{T}(k)=\ca{T}(k)\odot
(\psi_2^{-1}\odot \hat{\psi}_1)$ as desired. 

The bisection $\hat{\psi}_1^{-1}\odot \psi_2$ may be 
re-written, using 
\[\hat{\psi}_1^{-1}=\ca{T}\circ \psi_1^{-1}\circ (e_1^{-1}\circ
\ca{T}^*\circ e_{2})=\ca{T}\circ \psi_1^{-1}\circ \A(\psi_1)\circ \tau^*\circ \A(\psi_2)^{-1}\]
because of the commutativity of the diagram
\eqref{eq:commdia}. Thus, 
\[\hat{\psi}_1^{-1}(\A(\psi_2)(\mu))=\ca{T}\big(\psi_1(\tau^*\mu)\big)^{-1},\]
and \eqref{eq:anotherform} follows. 
\end{proof}

In the proof we used the following Lemma: 
\begin{lemma}\label{lem:twisttwist}
Suppose $M$ is a $K$-manifold with action $\A$. Let $\psi\in
\Gamma_\A(M,K)$ be a bisection relative to the action $\A$, and  
$\phi\in\Gamma_{\A^\psi}(M,K)$ a bisection relative to the twisted
action $\A^\psi$. Then the iterated twist $(\A^\psi)^\phi$ can be written as a
single twist, 
\[ (\A^\psi)^\phi=\A^{\psi\odot \A(\psi)^*\phi}.\]
\end{lemma}
\begin{proof}
By definition,
$\A^\psi(\phi)(x)=\A(\psi)\A(\phi(x))\A(\psi^{-1})(x)$. Hence
\[\A^\psi(\phi)=\A(\psi\odot \A(\psi)^*\phi\odot \psi^{-1}).\]
Using this formula we calculate, for all $k\in K$, 
\[ 
\begin{split}(\A^{\psi})^{\phi}(k)&=\A^{\psi}(\phi)\circ
\A^{\psi}(k)\circ \A^{\psi}(\phi)^{-1}\\
&=\A(\psi\odot \A(\psi)^*\phi)\circ \A(k)\circ
\A(\psi\odot \A(\psi)^*\phi)^{-1}\\
&=\A^{\psi\odot \A(\psi)^*\phi}(k).
\end{split}\]
\end{proof}

\subsection{Anti-Poisson involutions}\label{subsec:involutions}
An \emph{anti-Poisson involution} of a Poisson manifold $(M,\pi)$ is an
involutive diffeomorphism $s\in\on{Diff}(M)$ reversing the Poisson
structure, $s_*\pi=-\pi$. An anti-Poisson involution of a
Poisson Lie group $(K,\pi^K)$ is an anti-Poisson involution $s_K$ of
the underlying Poisson manifold which is also an automorphism of the
group $K$. In this case, $s_K$ canonically induces an anti-Poisson
involution of the dual Poisson Lie group $K^*$.

Suppose $K$ is a compact Lie group with standard
Poisson structure. Then any anti-linear involution of the Lie algebra $\g=\k^\C$ 
preserving the Iwasawa decomposition and the bilinear form $2\on{Im}B^\C$
defines an anti-Poisson involution $s_K$ of $K$. Let $s_{\k^*}$ be the induced involution of $\k^*$. 

\begin{lemma}
There exists a Ginzburg-Weinstein twist $\psi\in\Gamma(\k^*,K)$
which, in addition to the Properties from Theorem \ref{th:gw},
satisfies the equivariance property
\[ \psi\circ s_{\k^*}=s_K\circ\psi.\]
\end{lemma}
\begin{proof}
The Ginzburg-Weinstein twist constructed in the proof of Theorem 
\ref{th:gw} has the required equivariance property under involutions. 
Indeed, 
the forms $\sig_t,a_t$ on $\k^*$, hence also the Moser 1-form $b_t$, change
sign under $s_{\k^*}$ (see \cite{al:li}). 
It follows that the Moser vector field $v_t$ is
$s_{\k^*}$-invariant, while the function $\beta_t$ is equivariant, 
$\beta_t\circ s_{\k^*}=s_\k\circ \beta_t$. 
\end{proof}

\emph{Example.}
  If $K=\U(n)$, the complex conjugation operation $s_K(A)=\ol{A}$ is an
  anti-Poisson involution. The involution $s_{\k^*}$ is complex conjugation on 
  $\k^*\cong \Herm(n)$, and $s_{K^*}$ is 
  complex conjugation on upper triangular matrices or, equivalently,
  on the space $P=\Herm^+(n)$ of positive definite matrices.
Compatibility of a Ginzburg-Weinstein twist $\psi$ with these involutions 
just means $\psi(\ol{A})=\ol{\psi(A)}$. In particular, $\psi$  
restricts to a bisection $\on{Sym}(n)\to \SO(n)$.\\

The functoriality properties of Ginzburg-Weinstein maps generalize in
the obvious way to the presence of such involutions. Thus, suppose
$K_i,\ i=1,\ldots,n$ are compact Poisson Lie groups with standard Poisson
structure, and $s_{K_i}$ are anti-Poisson involutions of $K_i$ 
of the type discussed above.  Assume 
$\ca{T}_i\colon K_i\to K_{i+1},\ i=1\,\ldots,n-1$ are Poisson Lie
group homomorphisms with 
\[ \ca{T}_i\circ s_{K_i}=s_{K_{i+1}}\circ \ca{T}.\]
Then the Ginzburg-Weinstein twists $\psi_{i,t}\in\Gamma(\k_i^*,K_i)$
constructed in Theorem \ref{th:functorial} can be arranged to satisfy,
\[\psi_{i,t}\circ s_{\k_i^*}=s_{K_i}\circ \psi_{i,t}.\]
Indeed, the maps obtained in the proof of Theorem \ref{th:functorial}
automatically have this property, since all constructions are
compatible with the involutions.  It follows that all maps in the
commutative diagram \eqref{eq:commdia} intertwine the various
involutions. In particular, one obtains a commutative diagram for the
fixed point sets of the involutions.

\section{Gelfand-Zeitlin systems}

\subsection{Thimm actions}
The following construction of torus actions from non-Abelian group
actions appeared in Thimm's work \cite{th:in} on completely integrable
systems, and was later clarified by Guillemin-Sternberg in
\cite{gu:gc}. We will present the Thimm actions using the terminology
of bisections. Let $K$ be a compact Lie group, with maximal torus $T$,
and let $\k^*_{\on{reg}}\subset \k^*$ be the subset of regular
elements, that is, elements whose stabilizer is conjugate to $T$. Pick
a fundamental Weyl chamber $\t^*_+\subset \t^*$. Then
$\k^*_{\on{reg}}=K/T\times \on{int}(\t^*_+)$ as
$K$-manifolds. Restriction of equivariant bisections over
$\k^*_{\on{reg}}$ to $\on{int}(\t^*_+)$ defines a group isomorphism, 
\begin{equation}\label{eq:iso1}
\Gamma(\k^*_{\on{reg}},K)^K\xrightarrow{\cong} \Gamma(\on{int}(\t^*_+),T).
\end{equation}
\begin{lemma}
The isomorphism \eqref{eq:iso1} identifies the subgroups of Lagrangian
bisections: $\Gamma_0(\k^*_{\on{reg}},K)^K\cong \Gamma_0(\on{int}(\t^*_+),T)$. 
\end{lemma}
\begin{proof}
  If $\psi\in \Gamma(\k^*_{\on{reg}},K)^K$ is
  Lagrangian, then clearly so is its restriction to 
  $\on{int}(\t^*_+)$. For the converse, suppose $\psi$ restricts
  to a Lagrangian bisection over $\on{int}(\t^*_+)$.  For any
  $\xi\in\k$ we have
  $\iota(\xi)\l\mu,(\psi^{-1})^*\theta^L\r=\l\mu,\xi-\Ad_{\psi(\mu)}(\xi)\r=0$,
  since $\psi(\mu)\in K_\mu$. Hence also
\[ \iota(\xi)\d\l\mu,(\psi^{-1})^*\theta^L\r=
L(\xi)\l\mu,(\psi^{-1})^*\theta^L\r
-\d\iota(\xi)\l\mu,(\psi^{-1})^*\theta^L\r=0.\]
Since on the other hand the pull-back of
$\d\l\mu,(\psi^{-1})^*\theta^L\r$ to $\on{int}(\t^*_+)\subset
\k^*_{\on{reg}}$ is zero, this shows $\d
\l\mu,(\psi^{-1})^*\theta^L\r=0$. Thus $\psi$ is Lagrangian.
\end{proof}
Define a group homomorphism 
\begin{equation} 
\chi\colon T\to \Gamma_0(\k^*_{\on{reg}},K)^K
\end{equation}
by composing the map inverse to \eqref{eq:iso1} with the inclusion $T\to \Gamma_0(\on{int}(\t^*_+),T)$
as constant bisections. That is, $\chi(t)\colon \k^*_{\on{reg}}\to K$ is the unique
$K$-equivariant map with $\chi(t)(\mu)=t$ for $\mu\in\on{int}(\t^*_+)$.

Recall that by Lemma \ref{lem:center},
$\Gamma(\k^*_{\on{reg}},K)^K$ is the center of
$\Gamma(\k^*_{\on{reg}},K)$, and that its action on
$\k^*_{\on{reg}}$ is trivial. In particular, $\chi(t)$ acts trivially on 
$\k^*_{\on{reg}}$. Non-trivial actions are obtained by 
pulling $\chi(t)$ back under an equivariant map, 
$\Phi\colon M\to \k^*$. Thus let $M_0=\Phinv(\k^*_{\on{reg}})\subset
M$, and $\chi_M(t)=\Phi^*\chi(t)\in\Gamma(M_0,K)^K$. 
We define the \emph{Thimm action} of $t\in T$ by
\[ t\bullet x=\A(\chi_M(t))(x),\ \ x\in M_0\] 
By construction, the Thimm action commutes with the $K$-action, and 
the map $\Phi$ is Thimm-invariant:
\[ \Phi(t\bullet x)=t\bullet \Phi(x)=\Phi(x).\]
From now on, we will write $\chi(t)(\mu)\equiv \chi(t;\mu)$ and
similarly for $\chi_M$. 
\begin{lemma}\label{lem:easy1}
If $\psi\in \Gamma(M_0,K)$ is constant along the fibers of $\Phi$,
then $\psi$ commutes (under $\odot$) with all $\chi_M(t)$, and 
$(\psi\odot\chi_M(t))(x)=\psi(x)\chi_M(t;x)$.  
\end{lemma}
\begin{proof}
Since $\A(\chi_M(t))$ preserves the fibers of $\Phi$, the bisection 
$\psi$ satisfies $\A(\chi_M(t))^*\psi=\psi$. 
Hence, Lemma \ref{lem:easy} applies. 
\end{proof}
Thimm actions are naturally associated with Hamiltonian group actions. 
\begin{lemma}[Guillemin-Sternberg \cite{gu:gc}]
\label{lem:thimm}
Suppose $(M,\pi)$ is a Hamiltonian $K$-manifold, with moment map
$\Phi\colon M\to \k^*$. Then the Thimm $T$-action on $M_0$ 
is Hamiltonian, with moment map $q\circ \Phi\colon M\to \t^*$. 
Here $q\colon \k^*\to \t^*_+\subset \t^*$ is the unique $K$-invariant
map with $q(\mu)=\mu$ for $\mu\in\t^*_+$. 
\end{lemma}
Suppose now that 
\begin{equation}\label{eq:groups}
K_1\xrightarrow{\ca{T}_1} K_2
\xrightarrow{\ca{T}_2} \cdots \ra K_n
\end{equation}
 is a sequence of compact
Lie groups and homomorphisms, with differentials $\tau_i\colon \k_i\to
\k_{i+1}$. 
%
For $i<j$ we will write $\ca{T}_i^j=\ca{T}_{j-1}\circ\cdots\circ \ca{T}_i\colon
K_i\to K_j$, with differential $\tau_i^j\colon \k_i\to \k_j$. 
Take the maximal tori $T_i\subset K_i$ and positive Weyl
chambers $\t_{i,+}^*$ to be compatible, in the sense that for all
$i<n$, 
\[ \ca{T}_i(T_i)\subset T_{i+1},\ \ \tau_i^*(\t_{i+1,+}^*)\subset
\t_{i,+}^*.\]
Let $M$ be a $K_n$-manifold, and $\Phi_n\colon M\to \k_n^*$ an
equivariant map. Then each $K_i$ acts on $M$ via $\ca{T}_i^n$, and we
obtain a $K_i$-invariant map $\Phi_i=(\tau_i^n)^*\Phi_n\colon M\to
\k_i^*$. Let 
\[ M_0=\bigcap_{i=1}^n \Phi_i^{-1}(\k_{i,reg}^*),\]
and define $\chi_{i,M}\colon T_i\to \Gamma(M_0,K_n)$ by 
\[ \chi_{i,M}(t_i)=\ca{T}_i^n\circ \chi_i(t_i)\circ \Phi_i,\ \ t_i\in T_i\]
where $\chi_i(t_i)\in \Gamma(\k^*_{i,\on{reg}},K_i)$. 
\begin{lemma}
The images of the homomorphisms $\chi_{i,M}\colon T_i\to
\Gamma(M_0,K_n)$ all commute. Hence, they combine to define a group
homomorphism 
\[ \chi_M\colon T_n\times\cdots \times T_1\to  \Gamma(M_0,K_n).\]
One has the formula, 
\[ \chi_M(t_n,\ldots,t_1;x)=
\chi_{1,M}(t_1;x)\cdots \chi_{n,M}(t_n;x).
\]
\end{lemma}
\begin{proof}
Let $t_i\in T_i,\ t_j\in T_j$ where $i<j$. The bisection 
$\chi_j(t_j)\in \Gamma(\k^*_{j,\on{reg}},K_j)$ is $K_j$-equivariant, while
$\ca{T}_i^j\circ \chi_i(t_i)\circ (\tau_i^j)^*$ is constant along the fibers
of $(\tau_i^j)^*$. Hence, Lemma
\ref{lem:easy1} shows that the two bisections commute under $\odot$,
and that the product $(\ca{T}_i^j\circ \chi_i(t_i)\circ
(\tau_i^j)^*)\odot\chi_j(t_j)$ 
equals the pointwise product. It follows
that $\chi_{i,M}(t_i)$ and $\chi_{j,M}(t_j)$ commute and that the
product $\chi_{i,M}(t_i)\odot \chi_{j,M}(t_j)$ equals the pointwise
product. 
\end{proof}
We define the Thimm action of $t=(t_n,\ldots,t_1)\in T_n\times\cdots
\times T_1$ on $M_0$ by 
\[ t\bullet x=\A(\chi_M(t_n,\ldots,t_1))(x).\]
If $(M,\pi)$ is a Hamiltonian $K_n$-space, with moment map $\Phi_n$, 
then the Thimm action of $T_n\times \cdots
\times T_1$ on $M_0$ is Hamiltonian, with moment map 
\[ (q_n\circ \Phi_n,\ldots,q_1\circ \Phi_1)\colon M_0\to
\t^*_n\times\cdots \times \t^*_1.\]
Here $q_i\colon\k^*_{i}\to \t^*_{i,+}\subset \t^*_i$ are the
unique $K_i$-invariant maps with $q_i(\mu)=\mu$ for $\mu\in\t^*_+$.  
As a special case, the identity map $\Phi\colon \k_n^*\to \k_n^*$
gives rise to a Hamiltonian action of $T_{n-1}\times
\cdots \times T_1$ on 
\[(\k^*_n)_0=\bigcap_{i=1}^{n-1}
((\tau_i^n)^*)^{-1}(\k_{i,reg}^*).\]
(The torus $T_n$ is excluded, since its Thimm action is trivial.)

\subsection{Thimm actions for Poisson Lie groups}
Let $K$ be a compact Lie group with standard Poisson structure, and
$K^*_{\on{reg}}\subset K^*$ the subset 
of points whose stabilizer under the dressing action of $K$ 
has maximal rank. Since $e\colon \k^*_{\on{reg}}\to K^*_{\on{reg}}$ is 
a $K$-equivariant diffeomorphism, any $K$-equivariant map 
$\Psi\colon M\to K^*$ defines a Thimm $T$-action, via the composition 
$e^{-1}\circ \Psi$. Let $\psi\in\Gamma(\k^*,K)$ be a
Ginzburg-Weinstein twist, and $\gamma=e\circ \A(\psi)$. Parallel to Lemma \ref{lem:thimm} we have:
\begin{lemma}\label{lem:thimm1}
Suppose $M$ is a Poisson manifold, and 
$\Psi\colon M\to K^*$ is a moment map for a Poisson Lie group
action $\A\colon K\to \on{Diff}(M)$.  Then the Thimm $T$-action on 
$M_0$  is Hamiltonian, with moment map
\[ p\circ \Psi\colon M_0\to \t^*.\]
Here $p=q\circ e^{-1}\colon K^*\to \t^*$. If $\Psi=\gamma\circ \Phi$,
where $\Phi\colon M\to \k^*$ is a moment map for a Hamiltonian
$K$-action, then the Thimm actions defined by $\Phi$ and $\Psi$
coincide. 
\end{lemma}
\begin{proof}
  As shown in Proposition \ref{prop:moment}, $\Phi$ is the moment map
  for the twisted action $\A^{\Phi^*\psi}$ on $M$. Since
  $\A_{\k^*}(\psi)$ preserves orbits, $p=q\circ e^{-1}=q\circ
  \gamma^{-1}$. Thus, $p\circ \Psi=q\circ \Phi$ where
  $\Phi=\gamma^{-1}\circ \Psi$. Thus, Lemma \ref{lem:thimm} identifies
  $p\circ \Psi$ as the moment map for the Thimm $T$-action
  corresponding to $\Phi$ (relative to the \emph{twisted} action
  $\A^{\Phi^*\psi}$). Since the two $K$-actions are 
  conjugate under $\A(\Phi^*\psi)$, the same is true for the two Thimm
  $T$-actions.  But since $\chi_M(t)$ is $K$-equivariant, Lemma
  \ref{lem:easy1} shows $\Phi^*\psi\odot \chi_M(t)\circ
  \Phi^*\psi^{-1}=\chi_M(t)$. Hence, the two Thimm actions coincide.
\end{proof}

Suppose \eqref{eq:groups} is a sequence of homomorphisms of Poisson
Lie groups $K_1,\ldots,K_n$, equipped with the standard Poisson
structure.  Let $M$ be a $K_n$-manifold, let $\Psi_n\colon M\to K_n^*$
be an equivariant map, and let $\Psi_i\colon M\to K_i^*$ be the
composition of $\Psi_n$ with the map $(\ca{T}_i^n)^*\colon K_n^*\to
K_i^*$.  We then obtain commuting Thimm $T_i$-actions on
\[M_0=\bigcap_{i=1}^{n}\Psi_i^{-1}(K_{i,\on{reg}}^*).\]
If $(M,\pi)$ is a Poisson manifold, and $\Psi_n$ is the moment map for
a Poisson-Lie group action of $K_n$, then the Thimm $T_n \times \cdots \times T_1$-action on 
$M_0$  is Hamiltonian, with moment map
\[ (p_n\circ \Psi_n,\ldots,p_1\circ
\Psi_1)\colon M_0\to \t^*_n\times \cdots \times \t_1^*.\] Here
$p_i=q_i\circ e_i^{-1}$. In particular, we obtain a Hamiltonian
$T_{n-1}\times\cdots \times T_1$-action on
\[ (K_n^*)_0=\bigcap_{i=1}^{n-1} ((\ca{T}_i^n)^*)^{-1}(K_{i,\on{reg}}^*).\] 
By an inductive application of Theorem \ref{th:functorial}, it is possible
to choose Ginzburg-Weinstein twists $\psi_i\in\Gamma(\k_i^*,K_i)$, with
$\psi_i(0)=1$, which are \emph{compatible} in the sense the resulting diagram 
\[\label{eq:diagram0}
\begin{CD}
\k_n^* @>{\tau_{n-1}^*}>> \cdots @>{\tau_2^*}>> \k_2^* @>{\tau_1^*}>> \k_1^*\\
@VV{\gamma_n}V         @.      @VV{\gamma_2}V       @VV{\gamma_1}V \\
K_n^* @>>{{\ca{T}}_{n-1}^*}> \cdots @>>{{\ca{T}}_2^*}> K_2^* @>>{{\ca{T}}_1^*}> K_1^* 
\end{CD}
\]
with $\gamma_{i}=e_i\circ \A_i(\psi_i)$
commutes.

\begin{proposition}
For any choice of compatible Ginzburg-Weinstein twists
$\psi_i\in\Gamma(\k_i^*,K_i)$, 
the map $\gamma_n\colon \k_n^*\to K_n^*$ intertwines the Thimm 
$T_{n-1}\times\cdots \times T_1$-actions on $(\k_n^*)_0$ and
$(K_n^*)_0$, as well as their moment maps. 
The map $\psi_n$ has the following equivariance property under the 
Thimm action of $t=(t_{n-1},\ldots,t_1)\in T_{n-1}\times\cdots \times T_1$, 
\begin{equation}\label{eq:psinequiv}
\psi_n(t\bullet \mu)=\ti{\chi}(t;\mu)\,\psi_n(\mu)\,\chi(t;\mu)^{-1}.\end{equation}
Here 
\[ \chi(t;\mu)=\prod_{i=1}^{n-1} \ca{T}_i^n\big(\chi_i(t_i;\mu)\big),\ \ \
\ti{\chi}(t;\mu)=\prod_{i=1}^{n-1} \ca{T}_i^n\Big(\Ad_{\psi_i((\tau_i^n)^*\mu)}\chi_i(t_i;\mu)\Big).\]
\end{proposition}
\begin{proof}
For each $i<n$ we obtain commutative diagrams
\begin{equation}\label{eq:diagin}
\begin{CD}
\k_n^* @>{(\tau_i^n)^*}>> \k_{i}^* @>{q_i}>> \t^*_{i,+}\\
@VV{\gamma_n}V   @VV{\gamma_i}V       @VV{=}V \\
K_n^* @>>{(\ca{T}_i^n)^*}> K_i^* @>>{p_i}>   \t^*_{i,+}
\end{CD}
\end{equation}
It follows that the map $\gamma_n$ intertwines the moment maps for the
actions of $T_{n-1}\times\cdots \times T_1$, as well as the actions
themselves.  By Theorem \ref{th:equivariance}, the commutativity of
the diagram \eqref{eq:diagin} implies that the bisection
$\hat{\psi}_i^{-1}\odot \psi_n\in\Gamma(\k_n^*,K_n)$ is
$K_i$-equivariant, and that
\[(\hat{\psi}_i^{-1}\odot \psi_n)(\mu)
=\ca{T}_i^n\big(\psi_i((\tau_i^n)^*\mu)\big)^{-1}\,\psi_n(\mu).\]
The $K_i$-equivariance of the bisection $\psi=\hat{\psi}_i^{-1}\odot \psi_n$
implies the Thimm $T_i$-equivariance, 
\begin{equation}\label{eq:thimmeq} 
\psi(t_i\bullet\mu)=\Ad_{\ca{T}_i^n\circ
  \chi_i(t_i;\mu)}\psi(\mu).
\end{equation}
Using $(\tau_i^{n})^*(t_i\bullet\mu)=
(\tau_i^{n})^*\mu$, this yields
\[ \psi_n(t_i\bullet \mu)=
\ca{T}_i^n\Big(\Ad_{\psi_i((\tau_i^n)^*\mu)}\chi_i(t_i;\mu)\Big)\,\psi_n(\mu) \,\ca{T}_i^n\big(\chi_i(t_i;\mu)\big)^{-1},\]
proving \eqref{eq:psinequiv}.
\end{proof}

\emph{Remarks.}
\begin{enumerate}
\item
Throughout this discussion, we can assume that the functions $\psi_i$ 
take values in the semi-simple part $K_i^{ss}$. 
\item 
In the presence of anti-Poisson involutions $s_{K_i}$ (of the type
discussed in Section \ref{subsec:involutions}) with $s_{K_{i+1}}\circ
\ca{T}_i=\ca{T}_i\circ s_{K_{i}}$, one can assume that the maps
$\psi_i$ satisfy $s_{K_i}\circ \psi_i=\psi_i\circ s_{\k_i^*}$. Thus
$\gamma_n$ restricts to a diffeomorphism between the fixed point sets
of $s_{\k_n^*}$ and $s_{K_n^*}$, equivariant for the action of
$T_{n-1}'\times\cdots \times T_1'$, where $T_i'$ is the fixed point
set of the restriction of $s_{K_i}$ to $T_i$.
\end{enumerate}

\subsection{The $\U(n)$ Gelfand-Zeitlin system}\label{subsec:un}
Consider the sequence \eqref{eq:groups} for the special case
$K_i=\U(i)$, with the standard choice of maximal tori $T_i=T(i)$, and with $\ca{T}_i^j\colon U(i)\to \U(j)$
the inclusions as 
the upper left corner (extended by $1$'s along the
diagonal). Identifying $\mf{u}(i)^*\cong \Herm(i)$ as above, the 
standard choice of fundamental Weyl chamber consists of diagonal matrices with 
decreasing diagonal entries. The maps $(\tau_i^j)^*\colon
\mf{u}(j)^*\to \mf{u}(i)^*$ translate into the projection of 
a Hermitian $j\times j$-matrix onto the $i$th principal submatrix, and are clearly compatible with these
choices of $\t^*_{i,+}$. As shown by Guillemin-Sternberg
\cite{gu:gc}, the Thimm $T_{n-1}\times \cdots \times T_1$-action 
for the sequence of projections 
\[ \mf{u}(n)^*\to \cdots \to \mf{u}(2)^*\to \mf{u}(1)^*\]
defines a completely integrable system on $\mf{u}(n)^*$, and coincides
with the Gelfand-Zeitlin system described in Section \ref{sec:intro}. 

Let $U(i)$ carry the standard Poisson-Lie group structure
corresponding to these choices of $T_i,\ \t^*_{i,+}$ and the scalar
product $B_i(A',A)=-\on{tr}(A'A)$. The bracket on $\mf{u}(i)^*$ corresponds
to its identification with upper triangular matrices, with real
diagonal entries. The map 
\[ (\tau_i^j)^*\colon \mf{u}(j)^*\to \mf{u}(i)^*\] 
projects an upper triangular 
matrix onto the upper left $i\times i$ block, and is easily checked to
preserve Lie brackets. Hence, $\ca{T}_i^j$ are Poisson-Lie group
homomorphisms. The identification
\[ \U(i)^*\cong \on{Herm}^+(i),\]
takes the dressing action of $\U(i)$ to the action by conjugation. The maps 
\[(\ca{T}_i^j)^*\colon \U(j)^*\to \U(i)^*\] 
are again identified with projection to the upper left corner, both
under the identification with positive definite matrices, and under
the identification with the group upper triangular matrices with positive
diagonal.  The Thimm $T(n-1)\times \cdots \times T(1)$-action for the sequence of maps 
\[ \U(n)^*\to \cdots \to \U(2)^*\to \U(1)^*\]
is Flaschka-Ratiu's nonlinear Gelfand-Zeitlin system. Let
$\psi_i\colon \u(i)^*\to \SU(i)$ be compatible Ginzburg-Weinstein
twists, with $\psi_i(0)=1$ and $\psi_i(\ol{A})=\ol{\psi_i(A)}$, and
let $\gamma_i\colon \mf{u}(i)^*\to \U(i)^*$ be the corresponding
Ginzburg-Weinstein diffeomorphisms. Then $\psi_n\colon \u(n)^*\to
\SU(n)$ has the properties (i)-(iii) listed at the end of Section 
\ref{sec:linalg}. This finally completes the proof of Theorems 
\ref{th:real}, \ref{th:intertwines},
and \ref{th:phiext}. Furthermore, from the uniqueness properties 
of $\psi_n$ (Theorem \ref{th:phiext}), and since $\gamma_n$ is a 
Poisson map by construction, Theorem \ref{th:poisson} now comes for
free.\\

\emph{Remark.} While all the arguments in this paper were carried out
in the $C^\infty$-category, we could equally well have worked in the 
$C^\omega$-category of real-analytic maps. In particular, the distinguished 
2-form $\sig\in\Om^2(\k^*)$ from Section \ref{subsec:gw} is
real-analytic, by the explicit formula given in \cite{al:li}. 
It follows that the distinguished Ginzburg-Weinstein twist $\psi$ for 
$\U(n)$ is not only smooth, but is in fact real-analytic.\\

\subsection{Other classical groups}
We conclude with some remarks on Gelfand-Zeitlin systems for the
other classical groups. Consider first the special
orthogonal groups $\SO(n)$, with the standard choice of maximal tori. 
Guillemin-Sternberg's construction for the series of inclusions
\[ \SO(2)\to \SO(3)\to \cdots \]
produces a Gelfand-Zeitlin torus action over an open dense subset of
each Poisson manifold $\mf{so}(n)^*$. (Not to be confused with the
real locus of $\mf{u}(n)^*$, which does not carry a Poisson
structure.) A dimension count confirms that this defines a completely
integrable system. On the other hand, for the symplectic groups
the series of inclusions 
\[ \on{Sp}(1)\to \on{Sp}(2)\to \cdots \]
does not yield a completely integrable system, since the
Gelfand-Zeitlin torus does not have sufficiently large dimension. (By a
more sophisticated construction, Harada \cite{har:sym} was able to
obtain additional integrals of motion in this case.) Consider now the
standard Poisson structures on the groups $\SO(n)$ and $\on{Sp}(n)$.
Unfortunately, the inclusions $\SO(i)\to \SO(i+1)$ are \emph{not}
Poisson Lie group homomorphisms, essentially due to the fact that the
Dynkin diagram of $\SO(i)$ is not a subdiagram of that of $\SO(i+1)$.
However, the inclusions $\SO(i)\to \SO(i+2)$ are Poisson Lie group
homomorphisms, and so are the inclusions $\on{Sp}(i)\to \on{Sp}(i+1)$.
By the same discussion as for the unitary groups, one obtains
Ginzburg-Weinstein diffeomorphisms $\mf{so}(n)^*\to \SO(n)^*$ (resp.
$\mf{sp}(n)^*\to \on{Sp}(n)^*$) intertwining the resulting (partial)
Gelfand-Zeitlin systems. However, in contrast to the unitary groups,
there is no simple uniqueness statement in these cases.

\bibliographystyle{amsplain}

\def\polhk#1{\setbox0=\hbox{#1}{\ooalign{\hidewidth
  \lower1.5ex\hbox{`}\hidewidth\crcr\unhbox0}}} \def\cprime{$'$}
  \def\cprime{$'$} \def\polhk#1{\setbox0=\hbox{#1}{\ooalign{\hidewidth
  \lower1.5ex\hbox{`}\hidewidth\crcr\unhbox0}}} \def\cprime{$'$}
  \def\cprime{$'$}
\providecommand{\bysame}{\leavevmode\hbox to3em{\hrulefill}\thinspace}
\providecommand{\MR}{\relax\ifhmode\unskip\space\fi MR }
\providecommand{\MRhref}[2]{%
  \href{http://www.ams.org/mathscinet-getitem?mr=#1}{#2}
}
\providecommand{\href}[2]{#2}


\begin{thebibliography}{10}

\bibitem{al:po}
A.~Alekseev, \emph{On {P}oisson actions of compact {L}ie groups on symplectic
  manifolds}, J. Differential Geom. \textbf{45} (1997), no.~2, 241--256.

\bibitem{al:ka}
A.~Alekseev and E.~Meinrenken, \emph{Poisson geometry and the
  {K}ashiwara-{V}ergne conjecture}, C. R. Math. Acad. Sci. Paris \textbf{335}
  (2002), no.~9, 723--728.

\bibitem{al:li}
A.~Alekseev, E.~Meinrenken, and C.~Woodward, \emph{Linearization of {P}oisson
  actions and singular values of matrix products}, Ann. Inst. Fourier
  (Grenoble) \textbf{51} (2001), no.~6, 1691--1717.

\bibitem{bo:st}
P.~P. Boalch, \emph{Stokes matrices, {P}oisson {L}ie groups and {F}robenius
  manifolds}, Invent. Math. \textbf{146} (2001), no.~3, 479--506.

\bibitem{bu:ga}
H.~Bursztyn, \emph{On gauge transformations of poisson structures}, Quantum
  Field Theory and Noncommutative Geometry, Lecture Notes in Physics,
  Springer-Verlag.

\bibitem{bu:pi}
H.~Bursztyn and A. Weinstein, \emph{Picard groups in Poisson geometry},
Moscow Math. J. \textbf{4} (2004), 39-66. 

\bibitem{ch:gu}
V.~Chari and A.~Pressley, \emph{A guide to quantum groups}, Cambridge
  University Press, Cambridge, 1995, Corrected reprint of the 1994 original.

\bibitem{ca:ge}
A.~Cannas da~Silva and A.~Weinstein, \emph{Geometric models for noncommutative
  algebras}, American Mathematical Society, Providence, RI, 1999.

\bibitem{dr:qu}
V.~G. Drinfeld, \emph{Quantum groups}, Proceedings of the International
  Congress of Mathematicians, Vol. 1, 2 (Berkeley, Calif., 1986) (Providence,
  RI), Amer. Math. Soc., 1987, pp.~798--820.

\bibitem{du:on}
J.~J. Duistermaat and G.~J. Heckman, \emph{On the variation in the cohomology
  of the symplectic form of the reduced phase space}, Invent. Math. \textbf{69}
  (1982), 259--268.

\bibitem{en:co}
B.~Enriquez, P.~Etingof, and I.~Marshall, \emph{{Comparison of Poisson
  structures and Poisson-Lie dynamical r-matrices}}, Preprint, 2004.

\bibitem{fl:pc1}
H.~Flaschka and T.~Ratiu, \emph{A convexity theorem for {P}oisson actions of
  compact {L}ie groups}, IHES, Preprint 1995. Available at
  http://preprints.cern.ch.

\bibitem{fl:mo}
\bysame, \emph{A {M}orse-theoretic proof of {P}oisson {L}ie convexity},
  Integrable systems and foliations/Feuilletages et syst\`emes int\'egrables
  (Montpellier, 1995), Progr. Math., vol. 145, Birkh\"auser Boston, Boston, MA,
  1997, pp.~49--71.

\bibitem{gi:lp}
V.~L. Ginzburg and A.~Weinstein, \emph{Lie-{P}oisson structure on some
  {P}oisson {L}ie groups}, J. Amer. Math. Soc. \textbf{5} (1992), no.~2,
  445--453.

\bibitem{gu:gc}
V.~Guillemin and S.~Sternberg, \emph{The {G}elfand-{C}etlin system and
  quantization of the complex flag manifolds}, J. Funct. Anal \textbf{52}
  (1983), 106--128.

\bibitem{har:sym}
M.~Harada, \emph{{The symplectic geometry of the Gelfand-Cetlin-Molev basis
  for representations of Sp(2n,C)}}, Preprint, University of Toronto, 2004.
  {\tt arXiv:math.SG/0404485}.

\bibitem{ko:gel}
B.~Kostant and N.~Wallach, \emph{{Gelfand-Zeitlin theory from the perspective
  of classical mechanics. I}}, Preprint, {\tt arXiv:math.SG/0408342}.

\bibitem{lu:mo}
J.-H. Lu, \emph{Momentum mappings and reduction of {P}oisson actions},
  Symplectic geometry, groupoids, and integrable systems (Berkeley, CA, 1989),
  Springer, New York, 1991, pp.~209--226.

\bibitem{lu:po}
J.-H. Lu and A.~Weinstein, \emph{Poisson {L}ie groups, dressing
  transformations, and {B}ruhat decompositions}, J. Differential Geom.
  \textbf{31} (1990), no.~2, 501--526.

\bibitem{sev:po}
P.~{\v{S}}evera and A.~Weinstein, \emph{Poisson geometry with a 3-form
  background}, Progr. Theoret. Phys. Suppl. (2001), no.~144, 145--154,
  Noncommutative geometry and string theory (Yokohama, 2001).

\bibitem{th:in}
A.~Thimm, \emph{Integrable geodesic flows on homogeneous spaces}, Ergodic
  Theory Dynamical Systems \textbf{1} (1981), no.~4, 495--517 (1982).

\end{thebibliography}
\end{document}